%%%%%%%%%%%%%%%%%%%%%%%%%%%%%%%%%%%%%%%%%%%%%%%%%%%%%%%%%
%							%
%  Combinatorial Identities Related to Representations	%
%  of $U_q(\tilde{gl_2})$				%
%							%
%  Vitaly Tarasov					%
%							%
%  Preprint MPI 98-119, 8 pages				%
%							%
%  amstex.tex (ver. 2.1) and  amssym.tex  are required	%
%							%
%%%%%%%%%%%%%%%%%%%%%%%%%%%%%%%%%%%%%%%%%%%%%%%%%%%%%%%%%

%\mag 1200
\input amstex

\expandafter\ifx\csname beta.def\endcsname\relax \else\endinput\fi
\expandafter\edef\csname beta.def\endcsname{%
 \catcode`\noexpand\@=\the\catcode`\@\space}

\let\atbefore @

\catcode`\@=11

\overfullrule\z@
\hsize 6.25truein
\vsize 9.63truein

\let\@ft@\expandafter \let\@tb@f@\atbefore

\newif\ifMag
\ifnum\mag>1000 \Magtrue\fi
=\ifMag cmr8\else cmr9\fi

\newdimen\p@@ \p@@\p@
\def\m@ths@r{\ifnum\mathsurround=\z@\z@\else\maths@r\fi}
\def\maths@r{1.6\p@@} \def\mathsurzero{\def\maths@r{\z@}}

\mathsurround\maths@r
\font\Brm=cmr12 \font\Bbf=cmbx12 \font\Bit=cmti12 \font\ssf=cmss10
\font\Bsl=cmsl10 scaled 1200 \font\Bmmi=cmmi10 scaled 1200
\font\BBf=cmbx12 scaled 1200 \font\BMmi=cmmi10 scaled 1440

\def\atletter{\edef\atrestore{\catcode`\noexpand\@=\the\catcode`\@\space}
 \catcode`\@=11}

\newread\@ux \newwrite\@@x \newwrite\@@cd
\let\@np@@\input
\def\@np@t#1{\openin\@ux#1\relax\ifeof\@ux\else\closein\@ux\relax\@np@@ #1\fi}
\def\input#1 {\openin\@ux#1\relax\ifeof\@ux\wrs@x{No file #1}\else
 \closein\@ux\relax\@np@@ #1\fi}
\def\Input#1 {\relax} %% Do not remove the space after #1

\def\wr@@x#1{} \def\wrs@x{\immediate\write\sixt@@n}

\def\readldf{\@np@t{\jobname.ldf}}
\def\writeldf{\def\wr@@x{\immediate\write\@@x}
 \def\cl@selbl{\wr@@x{\string\Snodef{\the\Sno}}\wr@@x{\string\endinput}%
 \immediate\closeout\@@x} \immediate\openout\@@x\jobname.ldf}
\let\cl@selbl\relax

\def\tod@y{\ifcase\month\or
 January\or February\or March\or April\or May\or June\or July\or
 August\or September\or October\or November\or December\fi\space\,
\number\day,\space\,\number\year}

\newcount\c@time
\def\h@@r{hh}\def\m@n@te{mm}
\def\wh@tt@me{\c@time\time\divide\c@time 60\edef\h@@r{\number\c@time}%
 \multiply\c@time -60\advance\c@time\time\edef
 \m@n@te{\ifnum\c@time<10 0\fi\number\c@time}}
\def\t@me{\h@@r\/{\rm:}\m@n@te}  \let\whattime\wh@tt@me
\def\today{\tod@y\wr@@x{\string\todaydef{\tod@y}}}
\def\nowtime{\t@me{\let\/\ic@\wr@@x{\string\nowtimedef{\t@me}}}}
\def\todaydef#1{} \def\nowtimedef#1{}

\def\em#1{{\it #1\/}} \def\emph#1{{\sl #1\/}}

\def\fitem#1{\par\setbox\z@\hbox{#1}\hangindent\wd\z@
 \hglue-2\parindent\kern\wd\z@\indent\llap{#1}\ignore}

\def\itemflat#1{\par\setbox\z@\hbox{\rm #1\enspace}\hang\ifnum\wd\z@>\parindent
 \noindent\unhbox\z@\ignore\else\textindent{\rm#1}\fi}

\newcount\itemlet
\def\newbi{\itemlet 96} \newbi
\def\bitem{\gad\itemlet \par\hangindent1.5\parindent
 \hglue-.5\parindent\textindent{\rm\rlap{\char\the\itemlet}\hp{b})}}
\def\atem{\newbi\bitem}

\newcount\itemrm

\def\iitem{\gad\itemrm \par\hangindent1.5\parindent
 \hglue-.5\parindent\textindent{\rm\hp{v}\llap{\romannumeral\the\itemrm})}}

\def\center{\par\begingroup\leftskip\z@ plus \hsize \rightskip\leftskip
 \parindent\z@\parfillskip\z@skip \def\\{\unskip\break}}
\def\endcenter{\endgraf\endgroup}

\let\b@gr@@\begingroup \let\B@gr@@\begingroup
\def\b@gr@{\b@gr@@\let\b@gr@@\undefined}
\def\B@gr@{\B@gr@@\let\B@gr@@\undefined}

\def\@fn@xt#1#2#3{\let\@ch@r=#1\def\n@xt{\ifx\t@st@\@ch@r
 \def\n@@xt{#2}\else\def\n@@xt{#3}\fi\n@@xt}\futurelet\t@st@\n@xt}

\def\@fwd@@#1#2#3{\setbox\z@\hbox{#1}\ifdim\wd\z@>\z@#2\else#3\fi}
\def\s@twd@#1#2{\setbox\z@\hbox{#2}#1\wd\z@}

\def\r@st@re#1{\let#1\s@v@} \def\s@v@d@f{\let\s@v@}

\def\p@sk@p#1#2{\par\skip@#2\relax\ifdim\lastskip<\skip@\relax\removelastskip
 \ifnum#1=\z@\else\penalty#1\relax\fi\vskip\skip@
 \else\ifnum#1=\z@\else\penalty#1\relax\fi\fi}
\def\sk@@p#1{\par\skip@#1\relax\ifdim\lastskip<\skip@\relax\removelastskip
 \vskip\skip@\fi}

\newbox\p@b@ld
\def\poorbold#1{\setbox\p@b@ld\hbox{#1}\kern-.01em\copy\p@b@ld\kern-\wd\p@b@ld
 \kern.02em\copy\p@b@ld\kern-\wd\p@b@ld\kern-.012em\raise.02em\box\p@b@ld}

\ifx\plainfootnote\undefined \let\plainfootnote\footnote \fi

\let\s@v@\proclaim \let\proclaim\relax
\def\r@R@fs#1{\let#1\s@R@fs} \let\s@R@fs\Refs \let\Refs\relax
\def\r@endd@#1{\let#1\s@endd@} \let\s@endd@\enddocument
\let\bye\relax

\def\myR@fs{\@fn@xt[\m@R@f@\m@R@fs} \def\m@R@fs{\@fn@xt*\m@r@f@@\m@R@f@@}
\def\m@R@f@@{\m@R@f@[References]} \def\m@r@f@@*{\m@R@f@[]}

\def\Twelvepoint{\twelvepoint \let\Bbf\BBf \let\Bmmi\BMmi
\font\Brm=cmr12 scaled 1200 \font\Bit=cmti12 scaled 1200
\font\ssf=cmss10 scaled 1200 \font\Bsl=cmsl10 scaled 1440
\font\BBf=cmbx12 scaled 1440 \font\BMmi=cmmi10 scaled 1728}

\newif\ifamsppt

\newdimen\b@gsize

\def\p@@nt{.\kern.3em} \let\point\p@@nt

\let\proheadfont\bf \let\probodyfont\sl \let\demofont\it

\def\reffont#1{\def\r@ff@nt{#1}} \reffont\rm
\def\keyfont#1{\def\k@yf@nt{#1}} \keyfont\rm
\def\paperfont#1{\def\p@p@rf@nt{#1}} \paperfont\it
\def\bookfont#1{\def\b@@kf@nt{#1}} \bookfont\it
\def\volfont#1{\def\v@lf@nt{#1}} \volfont\bf
\def\issuefont#1{\def\iss@f@nt{#1}} \issuefont{no\p@@nt}

\newdimen\r@f@nd \newbox\r@f@b@x \newbox\adjb@x
\newbox\p@nct@ \newbox\k@yb@x \newcount\rcount
\newbox\b@b@x \newbox\p@p@rb@x \newbox\j@@rb@x \newbox\y@@rb@x
\newbox\v@lb@x \newbox\is@b@x \newbox\p@g@b@x \newif\ifp@g@ \newif\ifp@g@s
\newbox\inb@@kb@x \newbox\b@@kb@x \newbox\p@blb@x \newbox\p@bl@db@x
\newbox\ed@b@x \newif\ifed@ \newif\ifed@s \newif\if@fl@b \newif\if@fn@m
\newbox\p@p@nf@b@x \newbox\inf@b@x \newbox\b@@nf@b@x
\newtoks\@dd@p@n \newtoks\@ddt@ks

\newif\ifp@gen@

\@ft@\ifx\csname amsppt.sty\endcsname\relax

\headline={\hfil}
\footline={\ifp@gen@\ifnum\pageno=\z@\else\hfil\foliorm\folio\fi\else
 \ifnum\pageno=\z@\hfil\foliorm\folio\fi\fi\hfil\global\p@gen@true}
\parindent1pc

\font@\tensmc=cmcsc10
\font@\sevenex=cmex7
\font@\sevenit=cmti7
\font@\eightrm=cmr8
\font@\sixrm=cmr6
\font@\eighti=cmmi8 \skewchar\eighti='177
\font@\sixi=cmmi6 \skewchar\sixi='177
\font@\eightsy=cmsy8 \skewchar\eightsy='60
\font@\sixsy=cmsy6 \skewchar\sixsy='60
\font@\eightex=cmex8
\font@\eightbf=cmbx8
\font@\sixbf=cmbx6
\font@\eightit=cmti8
\font@\eightsl=cmsl8
\font@\eightsmc=cmcsc8
\font@\eighttt=cmtt8
\font@\ninerm=cmr9
\font@\ninei=cmmi9 \skewchar\ninei='177
\font@\ninesy=cmsy9 \skewchar\ninesy='60
\font@\nineex=cmex9
\font@\ninebf=cmbx9
\font@\nineit=cmti9
\font@\ninesl=cmsl9
\font@\ninesmc=cmcsc9
\font@\ninemsa=msam9
\font@\ninemsb=msbm9
\font@\nineeufm=eufm9
\font@\eightmsa=msam8
\font@\eightmsb=msbm8
\font@\eighteufm=eufm8
\font@\sixmsa=msam6
\font@\sixmsb=msbm6
\font@\sixeufm=eufm6

\loadmsam\loadmsbm\loadeufm
\input amssym.tex

\def\footnoterule{\kern-3\p@\hrule width5pc\kern 2.6\p@}
\def\m@k@foot#1{\insert\footins
 {\interlinepenalty\interfootnotelinepenalty
 \eightpoint\splittopskip\ht\strutbox\splitmaxdepth\dp\strutbox
 \floatingpenalty\@MM\leftskip\z@\rightskip\z@
 \spaceskip\z@\xspaceskip\z@
 \leavevmode\footstrut\ignore#1\unskip\lower\dp\strutbox
 \vbox to\dp\strutbox{}}}
\def\ftext#1{\m@k@foot{\vsk-.8>\nt #1}}
\def\pr@cl@@m#1{\p@sk@p{-100}\medskipamount
 \def\endproclaim{\endgroup\p@sk@p{55}\medskipamount}\begingroup
 \nt\ignore\proheadfont#1\unskip.\enspace\probodyfont\ignore}
\outer\def\proclaim{\pr@cl@@m} \s@v@d@f\proclaim \let\proclaim\relax
\def\demo#1{\sk@@p\medskipamount\nt{\ignore\demofont#1\unskip.}\enspace
 \ignore}
\def\enddemo{\sk@@p\medskipamount}

\def\cite#1{{\rm[#1]}} \let\nofrills\relax
 \def\Refs#1#2{\relax}

\def\big@#1#2{{\hbox{$\left#2\vcenter to#1\b@gsize{}%
 \right.\nulldelimiterspace\z@\m@th$}}}
\def\big{\big@\@ne}
\def\Big{\big@{1.5}}
\def\bigg{\big@\tw@}
\def\Bigg{\big@{2.5}}
\normallineskiplimit\p@

\def\tenpoint{\p@@\p@ \normallineskiplimit\p@@
 \mathsurround\m@ths@r \normalbaselineskip12\p@@
 \abovedisplayskip12\p@@ plus3\p@@ minus9\p@@
 \belowdisplayskip\abovedisplayskip
 \abovedisplayshortskip\z@ plus3\p@@
 \belowdisplayshortskip7\p@@ plus3\p@@ minus4\p@@
 \textonlyfont@\rm\tenrm \textonlyfont@\it\tenit
 \textonlyfont@\sl\tensl \textonlyfont@\bf\tenbf
 \textonlyfont@\smc\tensmc \textonlyfont@\tt\tentt
 \ifsyntax@ \def\big##1{{\hbox{$\left##1\right.$}}}%
  \let\Big\big \let\bigg\big \let\Bigg\big
 \else
  \textfont\z@\tenrm \scriptfont\z@\sevenrm \scriptscriptfont\z@\fiverm
  \textfont\@ne\teni \scriptfont\@ne\seveni \scriptscriptfont\@ne\fivei
  \textfont\tw@\tensy \scriptfont\tw@\sevensy \scriptscriptfont\tw@\fivesy
  \textfont\thr@@\tenex \scriptfont\thr@@\sevenex
	\scriptscriptfont\thr@@\sevenex
  \textfont\itfam\tenit \scriptfont\itfam\sevenit
	\scriptscriptfont\itfam\sevenit
  \textfont\bffam\tenbf \scriptfont\bffam\sevenbf
	\scriptscriptfont\bffam\fivebf
  \textfont\msafam\tenmsa \scriptfont\msafam\sevenmsa
	\scriptscriptfont\msafam\fivemsa
  \textfont\msbfam\tenmsb \scriptfont\msbfam\sevenmsb
	\scriptscriptfont\msbfam\fivemsb
  \textfont\eufmfam\teneufm \scriptfont\eufmfam\seveneufm
	\scriptscriptfont\eufmfam\fiveeufm
  \setbox\strutbox\hbox{\vrule height8.5\p@@ depth3.5\p@@ width\z@}%
  \setbox\strutbox@\hbox{\lower.5\normallineskiplimit\vbox{%
	\kern-\normallineskiplimit\copy\strutbox}}%
   \setbox\z@\vbox{\hbox{$($}\kern\z@}\b@gsize1.2\ht\z@
  \fi
  \normalbaselines\rm\dotsspace@1.5mu\ex@.2326ex\jot3\ex@}

\def\eightpoint{\p@@.8\p@ \normallineskiplimit\p@@
 \mathsurround\m@ths@r \normalbaselineskip10\p@
 \abovedisplayskip10\p@ plus2.4\p@ minus7.2\p@
 \belowdisplayskip\abovedisplayskip
 \abovedisplayshortskip\z@ plus3\p@@
 \belowdisplayshortskip7\p@@ plus3\p@@ minus4\p@@
 \textonlyfont@\rm\eightrm \textonlyfont@\it\eightit
 \textonlyfont@\sl\eightsl \textonlyfont@\bf\eightbf
 \textonlyfont@\smc\eightsmc \textonlyfont@\tt\eighttt
 \ifsyntax@\def\big##1{{\hbox{$\left##1\right.$}}}%
  \let\Big\big \let\bigg\big \let\Bigg\big
 \else
  \textfont\z@\eightrm \scriptfont\z@\sixrm \scriptscriptfont\z@\fiverm
  \textfont\@ne\eighti \scriptfont\@ne\sixi \scriptscriptfont\@ne\fivei
  \textfont\tw@\eightsy \scriptfont\tw@\sixsy \scriptscriptfont\tw@\fivesy
  \textfont\thr@@\eightex \scriptfont\thr@@\sevenex
	\scriptscriptfont\thr@@\sevenex
  \textfont\itfam\eightit \scriptfont\itfam\sevenit
	\scriptscriptfont\itfam\sevenit
  \textfont\bffam\eightbf \scriptfont\bffam\sixbf
	\scriptscriptfont\bffam\fivebf
  \textfont\msafam\eightmsa \scriptfont\msafam\sixmsa
	\scriptscriptfont\msafam\fivemsa
  \textfont\msbfam\eightmsb \scriptfont\msbfam\sixmsb
	\scriptscriptfont\msbfam\fivemsb
  \textfont\eufmfam\eighteufm \scriptfont\eufmfam\sixeufm
	\scriptscriptfont\eufmfam\fiveeufm
 \setbox\strutbox\hbox{\vrule height7\p@ depth3\p@ width\z@}%
 \setbox\strutbox@\hbox{\raise.5\normallineskiplimit\vbox{%
   \kern-\normallineskiplimit\copy\strutbox}}%
 \setbox\z@\vbox{\hbox{$($}\kern\z@}\b@gsize1.2\ht\z@
 \fi
 \normalbaselines\eightrm\dotsspace@1.5mu\ex@.2326ex\jot3\ex@}

\def\ninepoint{\p@@.9\p@ \normallineskiplimit\p@@
 \mathsurround\m@ths@r \normalbaselineskip11\p@
 \abovedisplayskip11\p@ plus2.7\p@ minus8.1\p@
 \belowdisplayskip\abovedisplayskip
 \abovedisplayshortskip\z@ plus3\p@@
 \belowdisplayshortskip7\p@@ plus3\p@@ minus4\p@@
 \textonlyfont@\rm\ninerm \textonlyfont@\it\nineit
 \textonlyfont@\sl\ninesl \textonlyfont@\bf\ninebf
 \textonlyfont@\smc\ninesmc \textonlyfont@\tt\ninett
 \ifsyntax@ \def\big##1{{\hbox{$\left##1\right.$}}}%
  \let\Big\big \let\bigg\big \let\Bigg\big
 \else
  \textfont\z@\ninerm \scriptfont\z@\sevenrm \scriptscriptfont\z@\fiverm
  \textfont\@ne\ninei \scriptfont\@ne\seveni \scriptscriptfont\@ne\fivei
  \textfont\tw@\ninesy \scriptfont\tw@\sevensy \scriptscriptfont\tw@\fivesy
  \textfont\thr@@\nineex \scriptfont\thr@@\sevenex
	\scriptscriptfont\thr@@\sevenex
  \textfont\itfam\nineit \scriptfont\itfam\sevenit
	\scriptscriptfont\itfam\sevenit
  \textfont\bffam\ninebf \scriptfont\bffam\sevenbf
	\scriptscriptfont\bffam\fivebf
  \textfont\msafam\ninemsa \scriptfont\msafam\sevenmsa
	\scriptscriptfont\msafam\fivemsa
  \textfont\msbfam\ninemsb \scriptfont\msbfam\sevenmsb
	\scriptscriptfont\msbfam\fivemsb
  \textfont\eufmfam\nineeufm \scriptfont\eufmfam\seveneufm
	\scriptscriptfont\eufmfam\fiveeufm
  \setbox\strutbox\hbox{\vrule height8.5\p@@ depth3.5\p@@ width\z@}%
  \setbox\strutbox@\hbox{\lower.5\normallineskiplimit\vbox{%
	\kern-\normallineskiplimit\copy\strutbox}}%
   \setbox\z@\vbox{\hbox{$($}\kern\z@}\b@gsize1.2\ht\z@
  \fi
  \normalbaselines\rm\dotsspace@1.5mu\ex@.2326ex\jot3\ex@}

\font@\twelverm=cmr10 scaled 1200
\font@\twelveit=cmti10 scaled 1200
\font@\twelvesl=cmsl10 scaled 1200
\font@\twelvebf=cmbx10 scaled 1200
\font@\twelvesmc=cmcsc10 scaled 1200
\font@\twelvett=cmtt10 scaled 1200
\font@\twelvei=cmmi10 scaled 1200 \skewchar\twelvei='177
\font@\twelvesy=cmsy10 scaled 1200 \skewchar\twelvesy='60
\font@\twelveex=cmex10 scaled 1200
\font@\twelvemsa=msam10 scaled 1200
\font@\twelvemsb=msbm10 scaled 1200
\font@\twelveeufm=eufm10 scaled 1200

\def\twelvepoint{\p@@1.2\p@ \normallineskiplimit\p@@
 \mathsurround\m@ths@r \normalbaselineskip12\p@@
 \abovedisplayskip12\p@@ plus3\p@@ minus9\p@@
 \belowdisplayskip\abovedisplayskip
 \abovedisplayshortskip\z@ plus3\p@@
 \belowdisplayshortskip7\p@@ plus3\p@@ minus4\p@@
 \textonlyfont@\rm\twelverm \textonlyfont@\it\twelveit
 \textonlyfont@\sl\twelvesl \textonlyfont@\bf\twelvebf
 \textonlyfont@\smc\twelvesmc \textonlyfont@\tt\twelvett
 \ifsyntax@ \def\big##1{{\hbox{$\left##1\right.$}}}%
  \let\Big\big \let\bigg\big \let\Bigg\big
 \else
  \textfont\z@\twelverm \scriptfont\z@\eightrm \scriptscriptfont\z@\sixrm
  \textfont\@ne\twelvei \scriptfont\@ne\eighti \scriptscriptfont\@ne\sixi
  \textfont\tw@\twelvesy \scriptfont\tw@\eightsy \scriptscriptfont\tw@\sixsy
  \textfont\thr@@\twelveex \scriptfont\thr@@\eightex
	\scriptscriptfont\thr@@\sevenex
  \textfont\itfam\twelveit \scriptfont\itfam\eightit
	\scriptscriptfont\itfam\sevenit
  \textfont\bffam\twelvebf \scriptfont\bffam\eightbf
	\scriptscriptfont\bffam\sixbf
  \textfont\msafam\twelvemsa \scriptfont\msafam\eightmsa
	\scriptscriptfont\msafam\sixmsa
  \textfont\msbfam\twelvemsb \scriptfont\msbfam\eightmsb
	\scriptscriptfont\msbfam\sixmsb
  \textfont\eufmfam\twelveeufm \scriptfont\eufmfam\eighteufm
	\scriptscriptfont\eufmfam\sixeufm
  \setbox\strutbox\hbox{\vrule height8.5\p@@ depth3.5\p@@ width\z@}%
  \setbox\strutbox@\hbox{\lower.5\normallineskiplimit\vbox{%
	\kern-\normallineskiplimit\copy\strutbox}}%
  \setbox\z@\vbox{\hbox{$($}\kern\z@}\b@gsize1.2\ht\z@
  \fi
  \normalbaselines\rm\dotsspace@1.5mu\ex@.2326ex\jot3\ex@}

\font@\twelvetrm=cmr10 at 12truept
\font@\twelvetit=cmti10 at 12truept
\font@\twelvetsl=cmsl10 at 12truept
\font@\twelvetbf=cmbx10 at 12truept
\font@\twelvetsmc=cmcsc10 at 12truept
\font@\twelvettt=cmtt10 at 12truept
\font@\twelveti=cmmi10 at 12truept \skewchar\twelveti='177
\font@\twelvetsy=cmsy10 at 12truept \skewchar\twelvetsy='60
\font@\twelvetex=cmex10 at 12truept
\font@\twelvetmsa=msam10 at 12truept
\font@\twelvetmsb=msbm10 at 12truept
\font@\twelveteufm=eufm10 at 12truept

\def\twelvetruepoint{\p@@1.2truept \normallineskiplimit\p@@
 \mathsurround\m@ths@r \normalbaselineskip12\p@@
 \abovedisplayskip12\p@@ plus3\p@@ minus9\p@@
 \belowdisplayskip\abovedisplayskip
 \abovedisplayshortskip\z@ plus3\p@@
 \belowdisplayshortskip7\p@@ plus3\p@@ minus4\p@@
 \textonlyfont@\rm\twelvetrm \textonlyfont@\it\twelvetit
 \textonlyfont@\sl\twelvetsl \textonlyfont@\bf\twelvetbf
 \textonlyfont@\smc\twelvetsmc \textonlyfont@\tt\twelvettt
 \ifsyntax@ \def\big##1{{\hbox{$\left##1\right.$}}}%
  \let\Big\big \let\bigg\big \let\Bigg\big
 \else
  \textfont\z@\twelvetrm \scriptfont\z@\eightrm \scriptscriptfont\z@\sixrm
  \textfont\@ne\twelveti \scriptfont\@ne\eighti \scriptscriptfont\@ne\sixi
  \textfont\tw@\twelvetsy \scriptfont\tw@\eightsy \scriptscriptfont\tw@\sixsy
  \textfont\thr@@\twelvetex \scriptfont\thr@@\eightex
	\scriptscriptfont\thr@@\sevenex
  \textfont\itfam\twelvetit \scriptfont\itfam\eightit
	\scriptscriptfont\itfam\sevenit
  \textfont\bffam\twelvetbf \scriptfont\bffam\eightbf
	\scriptscriptfont\bffam\sixbf
  \textfont\msafam\twelvetmsa \scriptfont\msafam\eightmsa
	\scriptscriptfont\msafam\sixmsa
  \textfont\msbfam\twelvetmsb \scriptfont\msbfam\eightmsb
	\scriptscriptfont\msbfam\sixmsb
  \textfont\eufmfam\twelveteufm \scriptfont\eufmfam\eighteufm
	\scriptscriptfont\eufmfam\sixeufm
  \setbox\strutbox\hbox{\vrule height8.5\p@@ depth3.5\p@@ width\z@}%
  \setbox\strutbox@\hbox{\lower.5\normallineskiplimit\vbox{%
	\kern-\normallineskiplimit\copy\strutbox}}%
  \setbox\z@\vbox{\hbox{$($}\kern\z@}\b@gsize1.2\ht\z@
  \fi
  \normalbaselines\rm\dotsspace@1.5mu\ex@.2326ex\jot3\ex@}

\font@\elevenrm=cmr10 scaled 1095
\font@\elevenit=cmti10 scaled 1095
\font@\elevensl=cmsl10 scaled 1095
\font@\elevenbf=cmbx10 scaled 1095
\font@\elevensmc=cmcsc10 scaled 1095
\font@\eleventt=cmtt10 scaled 1095
\font@\eleveni=cmmi10 scaled 1095 \skewchar\eleveni='177
\font@\elevensy=cmsy10 scaled 1095 \skewchar\elevensy='60
\font@\elevenex=cmex10 scaled 1095
\font@\elevenmsa=msam10 scaled 1095
\font@\elevenmsb=msbm10 scaled 1095
\font@\eleveneufm=eufm10 scaled 1095

\def\elevenpoint{\p@@1.1\p@ \normallineskiplimit\p@@
 \mathsurround\m@ths@r \normalbaselineskip12\p@@
 \abovedisplayskip12\p@@ plus3\p@@ minus9\p@@
 \belowdisplayskip\abovedisplayskip
 \abovedisplayshortskip\z@ plus3\p@@
 \belowdisplayshortskip7\p@@ plus3\p@@ minus4\p@@
 \textonlyfont@\rm\elevenrm \textonlyfont@\it\elevenit
 \textonlyfont@\sl\elevensl \textonlyfont@\bf\elevenbf
 \textonlyfont@\smc\elevensmc \textonlyfont@\tt\eleventt
 \ifsyntax@ \def\big##1{{\hbox{$\left##1\right.$}}}%
  \let\Big\big \let\bigg\big \let\Bigg\big
 \else
  \textfont\z@\elevenrm \scriptfont\z@\eightrm \scriptscriptfont\z@\sixrm
  \textfont\@ne\eleveni \scriptfont\@ne\eighti \scriptscriptfont\@ne\sixi
  \textfont\tw@\elevensy \scriptfont\tw@\eightsy \scriptscriptfont\tw@\sixsy
  \textfont\thr@@\elevenex \scriptfont\thr@@\eightex
	\scriptscriptfont\thr@@\sevenex
  \textfont\itfam\elevenit \scriptfont\itfam\eightit
	\scriptscriptfont\itfam\sevenit
  \textfont\bffam\elevenbf \scriptfont\bffam\eightbf
	\scriptscriptfont\bffam\sixbf
  \textfont\msafam\elevenmsa \scriptfont\msafam\eightmsa
	\scriptscriptfont\msafam\sixmsa
  \textfont\msbfam\elevenmsb \scriptfont\msbfam\eightmsb
	\scriptscriptfont\msbfam\sixmsb
  \textfont\eufmfam\eleveneufm \scriptfont\eufmfam\eighteufm
	\scriptscriptfont\eufmfam\sixeufm
  \setbox\strutbox\hbox{\vrule height8.5\p@@ depth3.5\p@@ width\z@}%
  \setbox\strutbox@\hbox{\lower.5\normallineskiplimit\vbox{%
	\kern-\normallineskiplimit\copy\strutbox}}%
  \setbox\z@\vbox{\hbox{$($}\kern\z@}\b@gsize1.2\ht\z@
  \fi
  \normalbaselines\rm\dotsspace@1.5mu\ex@.2326ex\jot3\ex@}

\def\m@R@f@[#1]{\mathsurzero{%\let\{\relax
 \s@ct{}{#1}}\wr@@c{\string\Refcd{#1}{\the\pageno}}\B@gr@
 \frenchspacing\rcount\z@\refkey{\k@yf@nt[##1]}\refno{\k@yf@nt[##1]}%
 \widest{AZ}\keyright\let\Key\key\let\refin\relax}
\def\widest#1{\s@twd@\r@f@nd{\r@fk@y{\k@yf@nt#1}\enspace}}
\def\widestno#1{\s@twd@\r@f@nd{\r@fn@{\k@yf@nt#1}\enspace}}
\def\widestlabel#1{\s@twd@\r@f@nd{\k@yf@nt#1\enspace}}
\def\refkey{\def\r@fk@y##1} \def\refno{\def\r@fn@##1}
\def\keyright{\def\r@fit@m{\hang\textindent}}
\def\keyflat{\def\r@fit@m##1{\setbox\z@\hbox{##1\enspace}\hang\noindent
 \ifnum\wd\z@<\parindent\indent\hglue-\wd\z@\fi\unhbox\z@}}

\def\R@fb@x{\global\setbox\r@f@b@x} \def\K@yb@x{\global\setbox\k@yb@x}
\def\ref{\par\b@gr@\r@ff@nt\R@fb@x\box\voidb@x\K@yb@x\box\voidb@x
 \@fn@mfalse\@fl@bfalse\b@g@nr@f}
\def\c@nc@t#1{\setbox\z@\lastbox
 \setbox\adjb@x\hbox{\unhbox\adjb@x\unhbox\z@\unskip\unskip\unpenalty#1}}
\def\adjust#1{\relax\ifmmode\penalty-\@M\null\hfil$\clubpenalty\z@
 \widowpenalty\z@\interlinepenalty\z@\offinterlineskip\endgraf
 \setbox\z@\lastbox\unskip\unpenalty\c@nc@t{#1}\nt$\hfil\penalty-\@M
 \else\endgraf\c@nc@t{#1}\nt\fi}
\def\adjustnext#1{\P@nct\hbox{#1}\ignore}
\def\adjustend#1{\def\@djp@{#1}\ignore}
\def\addtoks#1{\global\@ddt@ks{#1}\ignore}
\def\addnext#1{\global\@dd@p@n{#1}\ignore}

\def\cl@s@{\adjust{\@djp@}\endgraf\setbox\z@\lastbox
 \global\setbox\@ne\hbox{\unhbox\adjb@x\ifvoid\z@\else\unhbox\z@\unskip\unskip
 \unpenalty\fi}\egroup\ifnum\c@rr@nt=\k@yb@x\global\fi
 \setbox\c@rr@nt\hbox{\unhbox\@ne\box\p@nct@}\P@nct\null
 \the\@ddt@ks\global\@ddt@ks{}}
\def\@p@n#1{\def\c@rr@nt{#1}\setbox\c@rr@nt\vbox\bgroup\let\@djp@\relax
 \hsize\maxdimen\nt\the\@dd@p@n\global\@dd@p@n{}}
\def\b@g@nr@f{\bgroup\@p@n\z@}
\def\key{\cl@s@\ifvoid\k@yb@x\@p@n\k@yb@x\k@yf@nt\else\@p@n\z@\fi}
\def\label{\cl@s@\ifvoid\k@yb@x\global\@fl@btrue\@p@n\k@yb@x\k@yf@nt\else
 \@p@n\z@\fi}
\def\no{\cl@s@\ifvoid\k@yb@x\gad\rcount\global\@fn@mtrue
 \K@yb@x\hbox{\k@yf@nt\the\rcount}\fi\@p@n\z@}
\def\labelno{\cl@s@\ifvoid\k@yb@x\gad\rcount\@fl@btrue
 \@p@n\k@yb@x\k@yf@nt\the\rcount\else\@p@n\z@\fi}
\def\by{\cl@s@\@p@n\b@b@x} \def\paper{\cl@s@\@p@n\p@p@rb@x\p@p@rf@nt\ignore}
\def\jour{\cl@s@\@p@n\j@@rb@x} \def\yr{\cl@s@\@p@n\y@@rb@x}
\def\vol{\cl@s@\@p@n\v@lb@x\v@lf@nt\ignore}
\def\issue{\cl@s@\@p@n\is@b@x\iss@f@nt\ignore}
\def\page{\cl@s@\ifp@g@s\@p@n\z@\else\p@g@true\@p@n\p@g@b@x\fi}
\def\pages{\cl@s@\ifp@g@\@p@n\z@\else\p@g@strue\@p@n\p@g@b@x\fi}
\def\inbook{\cl@s@\@p@n\inb@@kb@x}
\def\book{\cl@s@\@p@n\b@@kb@x\b@@kf@nt\ignore}
\def\publ{\cl@s@\@p@n\p@blb@x} \def\publaddr{\cl@s@\@p@n\p@bl@db@x}
\def\ed{\cl@s@\ifed@s\@p@n\z@\else\ed@true\@p@n\ed@b@x\fi}
\def\eds{\cl@s@\ifed@\@p@n\z@\else\ed@strue\@p@n\ed@b@x\fi}
\def\info{\cl@s@\@p@n\inf@b@x} \def\paperinfo{\cl@s@\@p@n\p@p@nf@b@x}
\def\bookinfo{\cl@s@\@p@n\b@@nf@b@x} \let\finalinfo\info
\def\P@nct{\global\setbox\p@nct@} \def\nopunct{\P@nct\box\voidb@x}
\def\p@@@t#1#2{\ifvoid\p@nct@\else#1\unhbox\p@nct@#2\fi}
\def\sp@@{\penalty-50 \space\hskip\z@ plus.1em}
\def\c@mm@{\p@@@t,\sp@@} \def\sp@c@{\p@@@t\empty\sp@@}
\def\p@tb@x#1#2{\ifvoid#1\else#2\@nb@x#1\fi}
\def\@nb@x#1{\unhbox#1\P@nct\lastbox}
\def\endr@f@{\cl@s@\nopunct
 \R@fb@x\hbox{\unhbox\r@f@b@x \p@tb@x\b@b@x\empty
 \ifvoid\j@@rb@x\ifvoid\inb@@kb@x\ifvoid\p@p@rb@x\ifvoid\b@@kb@x
  \ifvoid\p@p@nf@b@x\ifvoid\b@@nf@b@x
  \p@tb@x\v@lb@x\c@mm@ \ifvoid\y@@rb@x\else\sp@c@(\@nb@x\y@@rb@x)\fi
  \p@tb@x\is@b@x\c@mm@ \p@tb@x\p@g@b@x\c@mm@ \p@tb@x\inf@b@x\c@mm@
  \else\p@tb@x \b@@nf@b@x\c@mm@ \p@tb@x\v@lb@x\c@mm@ \p@tb@x\is@b@x\sp@c@
  \ifvoid\ed@b@x\else\sp@c@(\@nb@x\ed@b@x,\space\ifed@ ed.\else eds.\fi)\fi
  \p@tb@x\p@blb@x\c@mm@ \p@tb@x\p@bl@db@x\c@mm@ \p@tb@x\y@@rb@x\c@mm@
  \p@tb@x\p@g@b@x{\c@mm@\ifp@g@ p\p@@nt\else pp\p@@nt\fi}%
  \p@tb@x\inf@b@x\c@mm@\fi
  \else \p@tb@x\p@p@nf@b@x\c@mm@ \p@tb@x\v@lb@x\c@mm@
  \ifvoid\y@@rb@x\else\sp@c@(\@nb@x\y@@rb@x)\fi
  \p@tb@x\is@b@x\c@mm@ \p@tb@x\p@g@b@x\c@mm@ \p@tb@x\inf@b@x\c@mm@\fi
  \else \p@tb@x\b@@kb@x\c@mm@
  \p@tb@x\b@@nf@b@x\c@mm@ \p@tb@x\p@blb@x\c@mm@
  \p@tb@x\p@bl@db@x\c@mm@ \p@tb@x\y@@rb@x\c@mm@
  \ifvoid\p@g@b@x\else\c@mm@\@nb@x\p@g@b@x p\fi \p@tb@x\inf@b@x\c@mm@ \fi
  \else \c@mm@\@nb@x\p@p@rb@x\ic@\p@tb@x\p@p@nf@b@x\c@mm@
  \p@tb@x\v@lb@x\sp@c@ \ifvoid\y@@rb@x\else\sp@c@(\@nb@x\y@@rb@x)\fi
  \p@tb@x\is@b@x\c@mm@ \p@tb@x\p@g@b@x\c@mm@\p@tb@x\inf@b@x\c@mm@\fi
  \else \p@tb@x\p@p@rb@x\c@mm@\ic@\p@tb@x\p@p@nf@b@x\c@mm@
  \c@mm@\@nb@x\inb@@kb@x \p@tb@x\b@@nf@b@x\c@mm@ \p@tb@x\v@lb@x\sp@c@
  \p@tb@x\is@b@x\sp@c@
  \ifvoid\ed@b@x\else\sp@c@(\@nb@x\ed@b@x,\space\ifed@ ed.\else eds.\fi)\fi
  \p@tb@x\p@blb@x\c@mm@ \p@tb@x\p@bl@db@x\c@mm@ \p@tb@x\y@@rb@x\c@mm@
  \p@tb@x\p@g@b@x{\c@mm@\ifp@g@ p\p@@nt\else pp\p@@nt\fi}%
  \p@tb@x\inf@b@x\c@mm@\fi
  \else\p@tb@x\p@p@rb@x\c@mm@\ic@\p@tb@x\p@p@nf@b@x\c@mm@\p@tb@x\j@@rb@x\c@mm@
  \p@tb@x\v@lb@x\sp@c@ \ifvoid\y@@rb@x\else\sp@c@(\@nb@x\y@@rb@x)\fi
  \p@tb@x\is@b@x\c@mm@ \p@tb@x\p@g@b@x\c@mm@ \p@tb@x\inf@b@x\c@mm@ \fi}}
\def\m@r@f#1#2{\endr@f@\ifvoid\p@nct@\else\R@fb@x\hbox{\unhbox\r@f@b@x
 #1\unhbox\p@nct@\penalty-200\enskip#2}\fi\egroup\b@g@nr@f}
\def\endref{\endr@f@\ifvoid\p@nct@\else\R@fb@x\hbox{\unhbox\r@f@b@x.}\fi
 \parindent\r@f@nd
 \r@fit@m{\ifvoid\k@yb@x\else\if@fn@m\r@fn@{\unhbox\k@yb@x}\else
 \if@fl@b\unhbox\k@yb@x\else\r@fk@y{\unhbox\k@yb@x}\fi\fi\fi}\unhbox\r@f@b@x
 \endgraf\egroup\endgroup}
\def\moreref{\m@r@f;\empty}
\def\transl{\m@r@f;{\unskip\space
 {\sl English translation\ic@}:\penalty-66 \space}}
\def\endRefs{\endgraf\goodbreak\endgroup}

\hyphenation{acad-e-my acad-e-mies af-ter-thought anom-aly anom-alies
an-ti-deriv-a-tive an-tin-o-my an-tin-o-mies apoth-e-o-ses
apoth-e-o-sis ap-pen-dix ar-che-typ-al as-sign-a-ble as-sist-ant-ship
as-ymp-tot-ic asyn-chro-nous at-trib-uted at-trib-ut-able bank-rupt
bank-rupt-cy bi-dif-fer-en-tial blue-print busier busiest
cat-a-stroph-ic cat-a-stroph-i-cally con-gress cross-hatched data-base
de-fin-i-tive de-riv-a-tive dis-trib-ute dri-ver dri-vers eco-nom-ics
econ-o-mist elit-ist equi-vari-ant ex-quis-ite ex-tra-or-di-nary
flow-chart for-mi-da-ble forth-right friv-o-lous ge-o-des-ic
ge-o-det-ic geo-met-ric griev-ance griev-ous griev-ous-ly
hexa-dec-i-mal ho-lo-no-my ho-mo-thetic ideals idio-syn-crasy
in-fin-ite-ly in-fin-i-tes-i-mal ir-rev-o-ca-ble key-stroke
lam-en-ta-ble light-weight mal-a-prop-ism man-u-script mar-gin-al
meta-bol-ic me-tab-o-lism meta-lan-guage me-trop-o-lis
met-ro-pol-i-tan mi-nut-est mol-e-cule mono-chrome mono-pole
mo-nop-oly mono-spline mo-not-o-nous mul-ti-fac-eted mul-ti-plic-able
non-euclid-ean non-iso-mor-phic non-smooth par-a-digm par-a-bol-ic
pa-rab-o-loid pa-ram-e-trize para-mount pen-ta-gon phe-nom-e-non
post-script pre-am-ble pro-ce-dur-al pro-hib-i-tive pro-hib-i-tive-ly
pseu-do-dif-fer-en-tial pseu-do-fi-nite pseu-do-nym qua-drat-ic
quad-ra-ture qua-si-smooth qua-si-sta-tion-ary qua-si-tri-an-gu-lar
quin-tes-sence quin-tes-sen-tial re-arrange-ment rec-tan-gle
ret-ri-bu-tion retro-fit retro-fit-ted right-eous right-eous-ness
ro-bot ro-bot-ics sched-ul-ing se-mes-ter semi-def-i-nite
semi-ho-mo-thet-ic set-up se-vere-ly side-step sov-er-eign spe-cious
spher-oid spher-oid-al star-tling star-tling-ly sta-tis-tics
sto-chas-tic straight-est strange-ness strat-a-gem strong-hold
sum-ma-ble symp-to-matic syn-chro-nous topo-graph-i-cal tra-vers-a-ble
tra-ver-sal tra-ver-sals treach-ery turn-around un-at-tached
un-err-ing-ly white-space wide-spread wing-spread wretch-ed
wretch-ed-ly Brown-ian Eng-lish Euler-ian Feb-ru-ary Gauss-ian
Grothen-dieck Hamil-ton-ian Her-mit-ian Jan-u-ary Japan-ese Kor-te-weg
Le-gendre Lip-schitz Lip-schitz-ian Mar-kov-ian Noe-ther-ian
No-vem-ber Rie-mann-ian Schwarz-schild Sep-tem-ber}

\def\leftheadtext#1{} \def\rightheadtext#1{}

\let\nopagenumber\p@gen@false \let\putpagenumber\p@gen@true
\let\pagefirst\nopagenumber \let\pagenext\putpagenumber

\else

\amsppttrue

\let\twelvepoint\relax \let\Twelvepoint\relax \let\putpagenumber\relax
\let\logo@\relax \let\pagefirst\firstpage@true \let\pagenext\firstpage@false
\def\nopagenumber{\let\f@li@ld\folio\def\folio{\global\let\folio\f@li@ld}}

\def\ftext#1{\footnotetext""{\vsk-.8>\nt #1}}

\def\m@R@f@[#1]{\Refs\nofrills{}\m@th\tenpoint
 {%\let\{\relax
 \s@ct{}{#1}}\wr@@c{\string\Refcd{#1}{\the\pageno}}
 \def\k@yf@##1{\hss[##1]\enspace} \let\keyformat\k@yf@
 \def\widest##1{\s@twd@\refindentwd{\tenpoint\k@yf@{##1}}}
 \let\Key\key \def\refin{\kern\refindentwd}}
\let\info\finalinfo \r@R@fs\Refs
\def\adjust#1{#1} \let\adjustend\relax
\let\adjustnext\adjust 

\fi

\outer\def\myRefs{\myR@fs} \r@st@re\proclaim
\def\bye{\par\vfill\supereject\cl@selbl\cl@secd\b@e} \r@endd@\b@e
\let\Cite\cite \let\Key\key \def\endpro{\par\endproclaim}
\let\d@c@\document \def\document{\d@c@\tenpoint}
\hyphenation{ortho-gon-al}

\newtoks\@@tp@t \@@tp@t\output
\output=\@ft@{\let\{\noexpand\the\@@tp@t}
\let\{\relax

\newif\ifVersion

\def\p@n@l#1{\ifnum#1=\z@\else\penalty#1\relax\fi}

\def\s@ct#1#2{\ifVersion
 \skip@\lastskip\ifdim\skip@<1.5\bls\vskip-\skip@\p@n@l{-200}\vsk.5>%
 \p@n@l{-200}\vsk.5>\p@n@l{-200}\vsk.5>\p@n@l{-200}\vsk-1.5>\else
 \p@n@l{-200}\fi\ifdim\skip@<.9\bls\vsk.9>\else
 \ifdim\skip@<1.5\bls\vskip\skip@\fi\fi
 \vtop{\twelvepoint\raggedright\s@cf@nt\vp1\vsk->\vskip.16ex
 \s@twd@\parindent{#1}%
 \ifdim\parindent>\z@\adv\parindent.5em\fi\hang\textindent{#1}#2\strut}
 \else
 \p@sk@p{-200}{.8\bls}\vtop{\s@cf@nt\s@twd@\parindent{#1}%
 \ifdim\parindent>\z@\adv\parindent.5em\fi\hang\textindent{#1}#2\strut}\fi
 \nointerlineskip\nobreak\vtop{\strut}\nobreak\vskip-.6\bls\nobreak}

\def\s@bs@ct#1#2{\ifVersion
 \skip@\lastskip\ifdim\skip@<1.5\bls\vskip-\skip@\p@n@l{-200}\vsk.5>%
 \p@n@l{-200}\vsk.5>\p@n@l{-200}\vsk.5>\p@n@l{-200}\vsk-1.5>\else
 \p@n@l{-200}\fi\ifdim\skip@<.9\bls\vsk.9>\else
 \ifdim\skip@<1.5\bls\vskip\skip@\fi\fi
 \vtop{\elevenpoint\raggedright\s@bf@nt\vp1\vsk->\vskip.16ex%
 \s@twd@\parindent{#1}\ifdim\parindent>\z@\adv\parindent.5em\fi
 \hang\textindent{#1}#2\strut}
 \else
 \p@sk@p{-200}{.6\bls}\vtop{\s@bf@nt\s@twd@\parindent{#1}%
 \ifdim\parindent>\z@\adv\parindent.5em\fi\hang\textindent{#1}#2\strut}\fi
 \nointerlineskip\nobreak\vtop{\strut}\nobreak\vskip-.8\bls\nobreak}

\def\gadv{\global\adv} \def\gad#1{\gadv#1\@ne} \def\gadneg#1{\gadv#1-\@ne}

\newcount\t@@n \t@@n=10 \newbox\testbox

\newcount\Sno \newcount\Lno \newcount\Fno

\def\pr@cl#1{\r@st@re\pr@c@\pr@c@{#1}\global\let\pr@c@\relax}

\def\tagg#1{\tag"\rlap{\rm(#1)}\kern.01\p@"}
\def\l@L#1{\l@bel{#1}L} \def\l@F#1{\l@bel{#1}F} \def\<#1>{\l@b@l{#1}F}
\def\Tag#1{\tag{\l@F{#1}}} \def\Tagg#1{\tagg{\l@F{#1}}}
\def\Rem{\demo{\sl Remark}} \def\Ex{\demo{\bf Example}}
\def\Pf#1.{\demo{Proof #1}} \def\epf{\qed\enddemo}
\def\Ap@x{Appendix}
\def\Appendix{\Sno=64 \t@@n\@ne \wr@@c{\string\Appencd}
 \def\sf@rm{\char\the\Sno} \def\sf@rm@{\Ap@x\space\sf@rm} \def\sf@rm@@{\Ap@x}
 \def\s@ct@n##1##2{\s@ct\empty{\setbox\z@\hbox{##1}\ifdim\wd\z@=\z@
 \if##2*\sf@rm@@\else\if##2.\sf@rm@@.\else##2\fi\fi\else
 \if##2*\sf@rm@\else\if##2.\sf@rm@.\else\sf@rm@.\enspace##2\fi\fi\fi}}}
\def\Appcd#1#2#3{\def\Ap@@{\hglue-\l@ftcd\Ap@x}\ifx\@ppl@ne\empty
 \def\l@@b{\@fwd@@{#1}{\space#1}{}}\if*#2\entcd{}{\Ap@@\l@@b}{#3}\else
 \if.#2\entcd{}{\Ap@@\l@@b.}{#3}\else\entcd{}{\Ap@@\l@@b.\enspace#2}{#3}\fi\fi
 \else\def\l@@b{\@fwd@@{#1}{\c@l@b{#1}}{}}\if*#2\entcd{\l@@b}{\Ap@x}{#3}\else
 \if.#2\entcd{\l@@b}{\Ap@x.}{#3}\else\entcd{\l@@b}{#2}{#3}\fi\fi\fi}

\let\s@ct@n\s@ct
\def\s@ct@@[#1]#2{\@ft@\xdef\csname @#1@S@\endcsname{\sf@rm}\wr@@x{}%
 \wr@@x{\string\labeldef{S}\space{\?#1@S?}\space{#1}}%
 {%\let\{\relax
 \s@ct@n{\sf@rm@}{#2}}\wr@@c{\string\Entcd{\?#1@S?}{#2}{\the\pageno}}}
\def\s@ct@#1{\wr@@x{}{%\let\{\relax
 \s@ct@n{\sf@rm@}{#1}}\wr@@c{\string\Entcd{\sf@rm}{#1}{\the\pageno}}}
\def\s@ct@e[#1]#2{\@ft@\xdef\csname @#1@S@\endcsname{\sf@rm}\wr@@x{}%
 \wr@@x{\string\labeldef{S}\space{\?#1@S?}\space{#1}}%
 {%\let\{\relax
 \s@ct@n\empty{#2}}\wr@@c{\string\Entcd{}{#2}{\the\pageno}}}
\def\s@cte#1{\wr@@x{}{%\let\{\relax
 \s@ct@n\empty{#1}}\wr@@c{\string\Entcd{}{#1}{\the\pageno}}}
\def\theSno#1#2{\dff\?#1@S?{#2}%
 \wr@@x{\string\labeldef{S}\space{#2}\space{#1}}\fi}

\newif\ifd@bn@\d@bn@true
\def\Section{\gad\Sno\ifd@bn@\Fno\z@\Lno\z@\fi\@fn@xt[\s@ct@@\s@ct@}
\def\section{\gad\Sno\ifd@bn@\Fno\z@\Lno\z@\fi\@fn@xt[\s@ct@e\s@cte}
\let\Sect\Section \let\sect\section
\def\subsection{\@fn@xt*\subs@ct@\subs@ct}
\def\subs@ct#1{{%\let\{\relax
 \s@bs@ct\empty{#1}}\wr@@c{\string\subcd{#1}{\the\pageno}}}
\def\subs@ct@*#1{\vsk->\vsk>{%\let\{\relax
 \s@bs@ct\empty{#1}}\wr@@c{\string\subcd{#1}{\the\pageno}}}
 \def\Snodef#1{\Sno #1}

\def\l@b@l#1#2{\def\n@@{\csname #2no\endcsname}%
 \if*#1\gad\n@@ \@ft@\xdef\csname @#1@#2@\endcsname{\l@f@rm}\else\def\t@st{#1}%
 \ifx\t@st\empty\gad\n@@ \@ft@\xdef\csname @#1@#2@\endcsname{\l@f@rm}%
 \else\@ft@\ifx\csname @#1@#2@mark\endcsname\relax\gad\n@@
 \@ft@\xdef\csname @#1@#2@\endcsname{\l@f@rm}%
 \@ft@\gdef\csname @#1@#2@mark\endcsname{}%
 \wr@@x{\string\labeldef{#2}\space{\?#1@#2?}\space\ifnum\n@@<10 \space\fi{#1}}%
 \fi\fi\fi}
\def\labeldef#1#2#3{\dff\?#3@#1?{#2}}
\def\Labeldef#1#2#3{\dff\?#3@#1?{#2}\@ft@\gdef\csname @#3@#1@mark\endcsname{}}
\def\Tagdef#1#2{\Labeldef{F}{#1}{#2}}

\def\l@bel#1#2{\l@b@l{#1}{#2}\?#1@#2?}

\newcount\c@cite
\def\?#1?{\csname @#1@\endcsname}
\def\[{\@fn@xt:\c@t@sect\c@t@}
\def\c@t@#1]{{\c@cite\z@\@fwd@@{\?#1@L?}{\adv\c@cite1}{}%
 \@fwd@@{\?#1@F?}{\adv\c@cite1}{}\@fwd@@{\?#1?}{\adv\c@cite1}{}%
 \relax\ifnum\c@cite=\z@{\bf ???}\wrs@x{No label [#1]}\else
 \ifnum\c@cite=1\let\@@PS\relax\let\@@@\relax\else\let\@@PS\underbar
 \def\@@@{{\rm<}}\fi\@@PS{\?#1?\@@@\?#1@L?\@@@\?#1@F?}\fi}}
\def\(#1){{\rm(\c@t@#1])}}
\def\c@t@s@ct#1{\@fwd@@{\?#1@S?}{\?#1@S?\relax}%
 {{\bf ???}\wrs@x{No section label {#1}}}}
\def\c@t@sect:#1]{\c@t@s@ct{#1}} \let\SNo\c@t@s@ct

\newdimen\l@ftcd \newdimen\r@ghtcd \let\nlc\relax

\def\d@tt@d{\leaders\hbox to 1em{\kern.1em.\hfil}\hfill}
\def\entcd#1#2#3{\item{\l@bcdf@nt#1}{\entcdf@nt#2}\alb\kern.9em\hbox{}%
 \kern-.9em\d@tt@d\kern-.36em{\p@g@cdf@nt#3}\kern-\r@ghtcd\hbox{}\par}
\def\Entcd#1#2#3{\def\l@@b{\@fwd@@{#1}{\c@l@b{#1}}{}}\vsk.2>%
 \entcd{\l@@b}{#2}{#3}}
\def\subcd#1#2{{\adv\leftskip.333em\entcd{}{\s@bcdf@nt#1}{#2}}}
\def\Refcd#1#2{\def\t@@st{#1}\ifx\t@@st\empty\ifx\r@fl@ne\empty\relax\else
 \R@fcd{\r@fl@ne}{#2}\fi\else\R@fcd{#1}{#2}\fi}
\def\R@fcd#1#2{\sk@@p{.6\bls}\entcd{}{\hglue-\l@ftcd\R@fcdf@nt #1}{#2}}
\def\Refline{\def\r@fl@ne} \def\Refempty{\let\r@fl@ne\empty}
\def\Appencd{\par\adv\leftskip-\l@ftcd\adv\rightskip-\r@ghtcd\@ppl@ne
 \adv\leftskip\l@ftcd\adv\rightskip\r@ghtcd\let\Entcd\Appcd}
\def\appline{\def\@ppl@ne} \def\Appempty{\let\@ppl@ne\empty}
\def\Appline#1{\def\@ppl@ne{\s@bs@ct{}{#1}}}
\def\leftcd#1{\adv\leftskip-\l@ftcd\s@twd@\l@ftcd{\c@l@b{#1}\enspace}
 \adv\leftskip\l@ftcd}
\def\rightcd#1{\adv\rightskip-\r@ghtcd\s@twd@\r@ghtcd{#1\enspace}
 \adv\rightskip\r@ghtcd}
\def\C@nt{Contents} \def\Ap@s{Appendices} \def\R@fcs{References}
\def\contents{\@fn@xt*\cont@@\cont@}
\def\cont@{\@fn@xt[\cnt@{\cnt@[\C@nt]}}
\def\cont@@*{\@fn@xt[\cnt@@{\cnt@@[\C@nt]}}
\def\cnt@[#1]{\c@nt@{M}{#1}{44}{\s@bs@ct{}{\@ppl@f@nt\Ap@s}}}
\def\cnt@@[#1]{\c@nt@{M}{#1}{44}{}}
\def\endco{\par\penalty-500\vsk>\vskip\z@\endgroup}
\def\readcd{\@np@t{\jobname.cd}}
\def\Cde{\@fn@xt*\Cde@@\Cde@}
\def\Cde@{\@fn@xt[\Cd@{\Cd@[\C@nt]}}
\def\Cde@@*{\@fn@xt[\Cd@@{\Cd@@[\C@nt]}}
\def\Cd@[#1]{\cnt@[#1]\readcd\endco}
\def\Cd@@[#1]{\cnt@@[#1]\readcd\endco}
\def\contlabeldef{\def\c@l@b}

\long\def\c@nt@#1#2#3#4{\s@twd@\l@ftcd{\c@l@b{#1}\enspace}
 \s@twd@\r@ghtcd{#3\enspace}\adv\r@ghtcd1.333em
 \def\@ppl@ne{#4}\def\r@fl@ne{\R@fcs}\s@ct{}{#2}\B@gr@\parindent\z@\let\nlc\nl
 \let\nl\relax\parskip.2\bls\adv\leftskip\l@ftcd\adv\rightskip\r@ghtcd}

\def\writecd{\immediate\openout\@@cd\jobname.cd \def\wr@@c{\write\@@cd}
 \def\cl@secd{\immediate\write\@@cd{\string\endinput}\immediate\closeout\@@cd}
 \def\closecd{\cl@secd\global\let\cl@secd\relax}}
\let\cl@secd\relax \def\wr@@c#1{} \let\closecd\relax

\def\dff{\@ft@\d@f} \def\d@f{\@ft@\def}
\def\edff{\@ft@\ed@f} \def\ed@f{\@ft@\edef}
\def\defi#1#2{\def#1{#2}\wr@@x{\string\def\string#1{#2}}}

\def\qed{\hbox{}\nobreak\hfill\nobreak{\m@th$\,\square$}}
\def\back#1 {\strut\kern-.33em #1\enspace\ignore} %% !!! a space after #1 !!!
\def\Text#1{\crcr\noalign{\alb\vsk>\normalbaselines\vsk->\vbox{\nt #1\strut}%
 \nobreak\nointerlineskip\vbox{\strut}\nobreak\vsk->\nobreak}}

\def\hcor#1{\advance\hoffset by #1}
\def\vcor#1{\advance\voffset by #1}
\let\bls\baselineskip \let\ignore\ignorespaces
\ifx\ic@\undefined \let\ic@\/\fi
\def\vsk#1>{\vskip#1\bls} \let\adv\advance
\def\vv#1>{\vadjust{\vsk#1>}\ignore}
\def\vvn#1>{\vadjust{\nobreak\vsk#1>\nobreak}\ignore}
\def\vvv#1>{\vskip\z@\vsk#1>\nt\ignore}
\def\vvgood{\vadjust{\penalty-500}} 
\def\Par{\vsk.5>} \def\setparindent{\edef\Parindent{\the\parindent}}
\def\Type{\vsk.5>\bgroup\parindent\z@\tt\rightskip\z@ plus1em minus1em%
 \spaceskip.3333em \xspaceskip.5em\relax}
\def\endType{\vsk.5>\egroup\nt} 

 \let\Tilde\widetilde \let\dollar\$ \let\ampersand\&
\let\sss\scriptscriptstyle  
\let\vp\vphantom \let\hp\hphantom \let\nt\noindent
\let\cline\centerline \let\lline\leftline \let\rline\rightline
\def\nn#1>{\noalign{\vskip#1\p@@}} \def\NN#1>{\openup#1\p@@}
\def\cnn#1>{\noalign{\vsk#1>}}
 
\let\Lim\lim \def\lim{\Lim\limits} \let\Sum\sum \def\sum{\Sum\limits}
 
\let\Prod\prod \def\prod{\Prod\limits} \let\Int\int \def\int{\Int\limits}

\def\tsum{\mathop{\tsize\Sum}\limits} 
\def\tprod{\mathop{\tsize\Prod}\limits}
\def\&{.\kern.1em} \def\>{{\!\;}} \def\]{{\!\!\;}} \def\){\>\]} \def\}{\]\]}
\def\nl{\leavevmode\hfill\break} \def\~{\leavevmode\@fn@xt~\m@n@s\@md@sh}
\def\m@n@s~{\raise.15ex\mbox{-}} \def\@md@sh{\raise.13ex\hbox{--}}
\let\procent\% \def\%#1{\ifmmode\mathop{#1}\limits\else\procent#1\fi}
\let\@ml@t\" \def\"#1{\ifmmode ^{(#1)}\else\@ml@t#1\fi}
\let\@c@t@\' \def\'#1{\ifmmode _{(#1)}\else\@c@t@#1\fi}
\let\colon\: \def\:{^{\vp|}}

\let\texspace\ \def\ {\ifmmode\alb\fi\texspace}
%%% Do not remove a space after \ !!!

\let\n@wp@ge\newpage \def\newpage{\endgraf\n@wp@ge}
\let\=\m@th \def\mbox#1{\hbox{\m@th$#1$}}
\def\mtext#1{\text{\m@th$#1$}} \def\^#1{\text{\m@th#1}}
\def\Line#1{\kern-.5\hsize\line{\m@th$\dsize#1$}\kern-.5\hsize}
\def\Lline#1{\kern-.5\hsize\lline{\m@th$\dsize#1$}\kern-.5\hsize}
\def\Cline#1{\kern-.5\hsize\cline{\m@th$\dsize#1$}\kern-.5\hsize}
\def\Rline#1{\kern-.5\hsize\rline{\m@th$\dsize#1$}\kern-.5\hsize}

\def\Ll@p#1{\llap{\m@th$#1$}} \def\Rl@p#1{\rlap{\m@th$#1$}}
 \def\Cl@p#1{\llap{\m@th$#1$\hss}}
\def\Llap#1{\mathchoice{\Ll@p{\dsize#1}}{\Ll@p{\tsize#1}}{\Ll@p{\ssize#1}}%
 {\Ll@p{\sss#1}}}
\def\Clap#1{\mathchoice{\Cl@p{\dsize#1}}{\Cl@p{\tsize#1}}{\Cl@p{\ssize#1}}%
 {\Cl@p{\sss#1}}}
\def\Rlap#1{\mathchoice{\Rl@p{\dsize#1}}{\Rl@p{\tsize#1}}{\Rl@p{\ssize#1}}%
 {\Rl@p{\sss#1}}}
 
\def\LRtph#1#2{\setbox\z@\hbox{#1}\dimen\z@\wd\z@\hbox{\hbox to\dimen\z@{#2}}}
\def\LRph#1#2{\LRtph{\m@th$#1$}{\m@th$#2$}}

 \def\RRph#1#2{\LRph{#1}{#2\hss}}

\def\Rph#1#2{\mathchoice{\RRph{\dsize#1}{\dsize#2}}{\RRph{\tsize#1}{\tsize#2}}
 {\RRph{\ssize#1}{\ssize#2}}{\RRph{\sss#1}{\sss#2}}}
\def\Lto#1{\setbox\z@\mbox{\tsize{#1}}%
 \mathrel{\mathop{\hbox to\wd\z@{\rightarrowfill}}\limits#1}}
\def\Lgets#1{\setbox\z@\mbox{\tsize{#1}}%
 \mathrel{\mathop{\hbox to\wd\z@{\leftarrowfill}}\limits#1}}
\def\vpb#1{{\vp{\big(}}^{\]#1}} \def\vpp#1{{\vp{\big]}}_{#1}}

\let\alb\allowbreak 
\def\ald{\noalign{\alb}} 

\let\o\circ \let\x\times \let\ox\otimes 
  \let\tabs\+
\let\le\leqslant \let\ge\geqslant
 \let\8\infty \let\*\star
\let\bra\langle \let\ket\rangle
 
\let\map\mapsto  
 
 \def\vert{\ |\ }

\let\lb\lbrace \let\rb\rbrace
 
\let\trileft\triangleleft \let\tright\triangleright

\def\lsym#1{#1\alb\ldots\relax#1\alb}
\def\lc{\lsym,}   \def\lox{\lsym\ox}
\def\llc{\,,\alb\ {\ldots\ ,}\alb\ }

\def\Res{\mathop{\roman{Res}\>}\limits}

\def\1{^{-1}} \def\_#1{_{\Rlap{#1}}}
\def\vst#1{{\lower1.9\p@@\mbox{\bigr|_{\raise.5\p@@\mbox{\ssize#1}}}}}
\def\vrp#1:#2>{{\vrule height#1 depth#2 width\z@}}
\def\vru#1>{\vrp#1:\z@>} \def\vrd#1>{\vrp\z@:#1>}
\def\qqq{\qquad\quad} 
\def\sscr#1{\raise.3ex\mbox{\sss#1}} \def\@@PS{\bold{OOPS!!!}}

\def\intcl{\mathop
 {\Rlap{\raise.3ex\mbox{\kern.12em\curvearrowleft}}\int}\limits}
\def\intcr{\mathop
 {\Rlap{\raise.3ex\mbox{\kern.24em\curvearrowright}}\int}\limits}

\def\pms{\raise.25ex\mbox{\ssize\pm}\>}
\def\mps{\raise.25ex\mbox{\ssize\mp}\>}

\let\al\alpha
\let\bt\beta
\let\gm\gamma  
\let\dl\delta \let\Dl\Delta 
\let\epe\epsilon \let\eps\varepsilon \let\epsilon\eps
\let\Xj\varXi

\let\tht\theta 
\let\thi\vartheta \let\Tho\varTheta

\let\ka\kappa
\let\la\lambda \let\La\Lambda

\let\si\sigma 
 
 \let\phi\varphi

\let\om\omega \let\Om\Omega 

\def\C{\Bbb C}

\def\Z{\Bbb Z}

\def\II{\Bbb I}

\def\TT{\Bbb T}

\def\Zp{\Z_{\ge 0}}

\def\difl/{differential} \def\dif/{difference}
\def\cf.{cf.\ \ignore} \def\Cf.{Cf.\ \ignore}
\def\egv/{eigenvector} \def\eva/{eigenvalue} \def\eq/{equation}
\def\lhs/{the left hand side} \def\rhs/{the right hand side}
\def\Lhs/{The left hand side} \def\Rhs/{The right hand side}
\def\gby/{generated by} \def\wrt/{with respect to} \def\st/{such that}
\def\resp/{respectively} \def\off/{offdiagonal} \def\wt/{weight}
\def\pol/{polynomial} \def\rat/{rational} \def\tri/{trigonometric}
\def\fn/{function} \def\var/{variable} \def\raf/{\rat/ \fn/}
\def\inv/{invariant} \def\hol/{holomorphic} \def\hof/{\hol/ \fn/}
\def\mer/{meromorphic} \def\mef/{\mer/ \fn/} \def\mult/{multiplicity}
\def\sym/{symmetric} \def\perm/{permutation} \def\fd/{finite-dimensional}
\def\rep/{representation} \def\irr/{irreducible} \def\irrep/{\irr/ \rep/}
\def\hom/{homomorphism} \def\aut/{automorphism} \def\iso/{isomorphism}
\def\lex/{lexicographical} \def\as/{asymptotic} \def\asex/{\as/ expansion}
\def\ndeg/{nondegenerate} \def\neib/{neighbourhood} \def\deq/{\dif/ \eq/}
\def\hw/{highest \wt/} \def\gv/{generating vector} \def\eqv/{equivalent}
\def\msd/{method of steepest descend} \def\pd/{pairwise distinct}
\def\wlg/{without loss of generality} \def\Wlg/{Without loss of generality}
\def\onedim/{one-dimensional} \def\qcl/{quasiclassical}
\def\hgeom/{hyper\-geometric} \def\hint/{\hgeom/ integral}
\def\hwm/{\hw/ module} \def\emod/{evaluation module} \def\Vmod/{Verma module}
\def\symg/{\sym/ group} \def\sol/{solution} \def\eval/{evaluation}
\def\anf/{analytic \fn/} \def\anco/{analytic continuation}
\def\qg/{quantum group} \def\qaff/{quantum affine algebra}

\def\Rm/{\^{$R$-}matrix} \def\Rms/{\^{$R$-}matrices} \def\YB/{Yang-Baxter \eq/}
\def\Ba/{Bethe ansatz} \def\Bv/{Bethe vector} \def\Bae/{\Ba/ \eq/}
\def\KZv/{Knizh\-nik-Zamo\-lod\-chi\-kov} \def\KZvB/{\KZv/-Bernard}
\def\KZ/{{\sl KZ\/}} \def\qKZ/{{\sl qKZ\/}}
\def\KZB/{{\sl KZB\/}} \def\qKZB/{{\sl qKZB\/}}
\def\qKZo/{\qKZ/ operator} \def\qKZc/{\qKZ/ connection}
\def\KZe/{\KZ/ \eq/} \def\qKZe/{\qKZ/ \eq/} \def\qKZBe/{\qKZB/ \eq/}

\def\h@ph{\discretionary{}{}{-}} \def\$#1$-{\,\^{$#1$}\h@ph}

\def\TFT/{Research Insitute for Theoretical Physics}
\def\HY/{University of Helsinki} \def\AoF/{the Academy of Finland}
\def\CNRS/{Supported in part by MAE\~MICECO\~CNRS Fellowship}
\def\LPT/{Laboratoire de Physique Th\'eorique ENSLAPP}
\def\ENSLyon/{\'Ecole Normale Sup\'erieure de Lyon}
\def\LPTaddr/{46, All\'ee d'Italie, 69364 Lyon Cedex 07, France}
\def\enslapp/{URA 14\~36 du CNRS, associ\'ee \`a l'E.N.S.\ de Lyon,
au LAPP d'Annecy et \`a l'Universit\`e de Savoie}
\def\ensemail/{vtarasov\@ enslapp.ens-lyon.fr}
\def\DMS/{Department of Mathematics, Faculty of Science}
\def\DMO/{\DMS/, Osaka University}
\def\DMOaddr/{Toyonaka, Osaka 560, Japan}
\def\dmoemail/{vt\@ math.sci.osaka-u.ac.jp}
\def\SPb/{St\&Peters\-burg}
\def\home/{\SPb/ Branch of Steklov Mathematical Institute}
\def\homeaddr/{Fontanka 27, \SPb/ \,191011, Russia}
\def\homemail/{vt\@ pdmi.ras.ru}
\def\absence/{On leave of absence from \home/}
\def\UNC/{Department of Mathematics, University of North Carolina}
\def\ChH/{Chapel Hill}
\def\UNCaddr/{\ChH/, NC 27599, USA} \def\avemail/{av\@ math.unc.edu}
\def\grant/{NSF grant DMS\~9501290}	%%% Felder's grant no. 9400841
\def\Grant/{Supported in part by \grant/}

\def\Aomoto/{K\&Aomoto}
\def\Dri/{V\]\&G\&Drin\-feld}
\def\Fadd/{L\&D\&Fad\-deev}
\def\Feld/{G\&Felder}
\def\Fre/{I\&B\&Fren\-kel}
\def\Gustaf/{R\&A\&Gustafson}
\def\Kazh/{D\&Kazhdan} \def\Kir/{A\&N\&Kiril\-lov}
\def\Kor/{V\]\&E\&Kore\-pin}
\def\Lusz/{G\&Lusztig}
\def\MN/{M\&Naza\-rov}
\def\Resh/{N\&Reshe\-ti\-khin} \def\Reshy/{N\&\]Yu\&Reshe\-ti\-khin}
\def\SchV/{V\]\&\]V\]\&Schecht\-man} \def\Sch/{V\]\&Schecht\-man}
\def\Skl/{E\&K\&Sklya\-nin}
\def\Smirn/{F\]\&Smirnov} \def\Smirnov/{F\]\&A\&Smirnov}
\def\Takh/{L\&A\&Takh\-tajan}
\def\VT/{V\]\&Ta\-ra\-sov} \def\VoT/{V\]\&O\&Ta\-ra\-sov}
\def\Varch/{A\&\]Var\-chenko} \def\Varn/{A\&N\&\]Var\-chenko}

\def\AMS/{Amer.\ Math.\ Society}
\def\CMP/{Comm.\ Math.\ Phys.{}}
\def\DMJ/{Duke.\ Math.\ J.{}}
\def\Inv/{Invent.\ Math.{}} %% Inventiones Mathematicae
\def\IMRN/{Int.\ Math.\ Res.\ Notices}
\def\JPA/{J.\ Phys.\ A{}}
\def\JSM/{J.\ Soviet\ Math.{}}
\def\LMP/{Lett.\ Math.\ Phys.{}}
\def\LMJ/{Leningrad Math.\ J.{}}
\def\LpMJ/{\SPb/ Math.\ J.{}}
\def\SIAM/{SIAM J.\ Math.\ Anal.{}}
\def\SMNS/{Selecta Math., New Series}
\def\TMP/{Theor.\ Math.\ Phys.{}}
\def\ZNS/{Zap.\ nauch.\ semin. LOMI}

\def\ASMP/{Advanced Series in Math.\ Phys.{}}

\def\AMSa/{AMS \publaddr Providence}
\def\Birk/{Birkh\"auser}
\def\CUP/{Cambridge University Press} \def\CUPa/{\CUP/ \publaddr Cambridge}
\def\Spri/{Springer-Verlag} \def\Spria/{\Spri/ \publaddr Berlin}
\def\WS/{World Scientific} \def\WSa/{\WS/ \publaddr Singapore}

\newbox\lefthbox \newbox\righthbox

\let\sectsep. \let\labelsep. \let\contsep. \let\labelspace\relax
\let\sectpre\relax \let\contpre\relax
\def\sf@rm{\the\Sno} \def\sf@rm@{\sectpre\sf@rm\sectsep}
\def\c@l@b#1{\contpre#1\contsep}
\def\l@f@rm{\ifd@bn@\sf@rm\labelsep\fi\labelspace\the\n@@}

\def\sectformdef{\def\sf@rm}

\let\DoubleNum\d@bn@true \let\SingleNum\d@bn@false

\def\NoNewNum{\let\writeldf\relax\def\l@b@l##1##2{\if*##1%
 \@ft@\xdef\csname @##1@##2@\endcsname{\mbox{*{*}*}}\fi}}
\def\NoNewTime{\def\todaydef##1{\def\today{##1}}
 \def\nowtimedef##1{\def\nowtime{##1}}}
\def\NoInput{\let\Input\input\let\writeldf\relax}
\def\Fixed{\NoNewTime\NoInput}

\def\sectfont#1{\def\s@cf@nt{#1}} \sectfont\bf
\def\subsectfont#1{\def\s@bf@nt{#1}} \subsectfont\it
\def\Entcdfont#1{\def\entcdf@nt{#1}} \Entcdfont\relax
\def\labelcdfont#1{\def\l@bcdf@nt{#1}} \labelcdfont\relax
\def\pagecdfont#1{\def\p@g@cdf@nt{#1}} \pagecdfont\relax
\def\subcdfont#1{\def\s@bcdf@nt{#1}} \subcdfont\it
\def\applefont#1{\def\@ppl@f@nt{#1}} \applefont\bf
\def\Refcdfont#1{\def\R@fcdf@nt{#1}} \Refcdfont\bf

\tenpoint

\csname beta.def\endcsname

\Fixed

\expandafter\ifx\csname ident.def\endcsname\relax \else\endinput\fi
\expandafter\edef\csname ident.def\endcsname{%
 \catcode`\noexpand\@=\the\catcode`\@\space}

\catcode`\@=11

\font@\Beufm=eufm10 scaled 1440

\newif\ifMPIM

\def\Th#1{\pr@cl{Theorem \l@L{#1}}\ignore}
\def\Lm#1{\pr@cl{Lemma \l@L{#1}}\ignore}
\def\Cr#1{\pr@cl{Corollary \l@L{#1}}\ignore}
\def\Df#1{\pr@cl{Definition \l@L{#1}}\ignore}
\def\Cj#1{\pr@cl{Conjecture \l@L{#1}}\ignore}
\def\Prop#1{\pr@cl{Proposition \l@L{#1}}\ignore}

\def\EE{\Bbb E}
\def\Sb{\bold S}
\def\Ec{\Cal E}
\def\Pc{\Cal P}
\def\lg{\frak l}

\def\vti{\tilde v}

\def\Cl{\C^{\>\ell}}
\def\qH{q^{-H}} \def\q{q\1}
\def\Ae{A^{\sscr{e\]l\}l}}}

\def\Sl{\Sb_\ell}
\def\Pln{\Pc_{\ell,n}} \def\Plij{\Pc_\ell^{i,j}}
\def\susi{\sum_{\,\si\in\Sl\!}}
\def\sula{\sum_{\,\la\in\Pln\!}}
\def\sumu{\sum_{\,\mu\in\Pln\!}}

\def\l@inf{\lower.21ex\mbox{\ssize\8}} \def\9{_{\kern-.02em\l@inf}}
\def\ts{t^{\sss\star}} \def\dtt{(dt/t)^\ell}
\def\Sf{_{\]S}\:} \def\Of{_\Om\:}

 \def\tpron{\tprod_{m=1}^n}
\def\sun{\sum_{m=1}^n} \def\pron{\prod_{m=1}^n} \def\prmn{\prod_{m=1}^{n-1}}
\def\pral{\prod_{a=1}^\ell} 
\def\prab{\prod_{1\le a<b\le\ell}}
\def\prjm{\prod_{1\le j<m}} \def\tprjm{\!\tprod_{1\le j<m}\!}
\def\pros{\prod_{s=0}^{\ell-1}} \def\prsl{\prod_{s=1-\ell}^{\ell-1}}

\def\mn{m=1\lc n} \def\aell{a=1\lc\ell}
\def\xn{x_1\lc x_n} \def\yn{y_1\lc y_n} \def\zn{z_1\lc z_n}
\def\tell{t_1\lc t_\ell}  
\def\xt{x\mathbin{\smash\tright\]}} \def\yt{y\mathbin{\]\smash\trileft}}
\def\xRes{\mathop{\xt\]\italic{Res}\)}\limits}
\def\yRes{\mathop{\yt\]\italic{Res}\)}\limits}

\def\vox{v_1\lox v_n} \def\voxx{v_n\lox v_1}
\def\Voz{V_1(z_1)\lox V_n(z_n)} \def\Vozz{V_n(z_n)\lox V_1(z_1)}
\def\Fv{F\vpb{\,\om_{1,\la}}v_1\lox F\vpb{\,\om_{n,\la}}v_n}

\def\gsl{\frak{sl}_2} \def\gl{\frak{gl}_2}
\def\Uq{U_q(\gsl)} \def\Uqh{U_q(\Tilde{\gl})}
\def\Eqg{E_{\rho,\gm}(\gsl)}

\def\indet/{indeterminate}
\def\qlo/{quantum loop algebra} \def\eqg/{elliptic \qg/}
\def\Umod/{\$\Uq$-module} \def\Uhmod/{\$\Uqh$-module}
\def\Emod/{\$\Eqg$-module} \def\p-{\^{$\,p\)$}-}

\let\goodbm\relax  \let\mmgood\relax 
   
 \def\vvm#1>{\relax} 

\ifMag\let\goodbm\goodbreak  \let\mmgood\vvgood
 \let\goodbreak\relax  \let\vvgood\relax 
  \let\vvm\vv \fi

\csname ident.def\endcsname

\labeldef{S} {1} {Pi}
\labeldef{F} {1\labelsep \labelspace 1}  {P}
\labeldef{F} {1\labelsep \labelspace 2}  {P'}
\labeldef{L} {1\labelsep \labelspace 1}  {id}
\labeldef{F} {1\labelsep \labelspace 3}  {id1}

\labeldef{S} {2} {Pf}
\labeldef{L} {2\labelsep \labelspace 1}  {PP}
\labeldef{F} {2\labelsep \labelspace 1}  {N}
\labeldef{L} {2\labelsep \labelspace 2}  {Pxy}
\labeldef{F} {2\labelsep \labelspace 2}  {xRes}
\labeldef{L} {2\labelsep \labelspace 3}  {ResI}
\labeldef{F} {2\labelsep \labelspace 3}  {<>}
\labeldef{F} {2\labelsep \labelspace 4}  {PQ}
\labeldef{F} {2\labelsep \labelspace 5}  {MN}
\labeldef{F} {2\labelsep \labelspace 6}  {id2}
\labeldef{F} {2\labelsep \labelspace 7}  {detQ}
\labeldef{F} {2\labelsep \labelspace 8}  {detA}

\labeldef{S} {3} {Uq}
\labeldef{F} {3\labelsep \labelspace 1}  {Ughr}
\labeldef{L} {3\labelsep \labelspace 1}  {BC}
\labeldef{L} {3\labelsep \labelspace 2}  {BCC}
\labeldef{F} {3\labelsep \labelspace 2}  {BC1}
\labeldef{F} {3\labelsep \labelspace 3}  {BC2}
\labeldef{L} {3\labelsep \labelspace 3}  {KBI}

\labeldef{S} {4} {Ei}
\labeldef{F} {4\labelsep \labelspace 1}  {T}
\labeldef{F} {4\labelsep \labelspace 2}  {T'}
\labeldef{L} {4\labelsep \labelspace 1}  {idp}
\labeldef{F} {4\labelsep \labelspace 3}  {idp1}
\labeldef{F} {4\labelsep \labelspace 4}  {idp2}
\labeldef{L} {4\labelsep \labelspace 2}  {Res}
\labeldef{L} {4\labelsep \labelspace 3}  {XX}
\labeldef{F} {4\labelsep \labelspace 5}  {D}
\labeldef{F} {4\labelsep \labelspace 6}  {thi}
\labeldef{F} {4\labelsep \labelspace 7}  {Tho}
\labeldef{F} {4\labelsep \labelspace 8}  {detT}
\labeldef{F} {4\labelsep \labelspace 9}  {detAe}
\labeldef{L} {4\labelsep \labelspace 4}  {XT}

\Versiontrue

%\ifMag\MPIMtrue\fi

\ifMPIM
\hfuzz 15pt
\center
\vglue0pt
{\Twelvepoint \bls1.44\bls\Brm
Combinatorial Identities Related to Representations of
$\hbox{\Bmmi U}_{\tsize q}\hbox{(\mbox{\Tilde{\hbox{\Beufm gl}_{\tsize2}}})}$
\vsk1.07>
\smc \VT/}
\vsk1.7>
{\sl October \,16, 1998}
\vsk>
\endcenter
\ftext{\=\tenpoint\sl E-mail\/{\rm:} \homemail/}
\else
\headline={\hfil\ifnum\pageno=1\hbox{\tenpoint MPI \,98\~~119}\fi}
\center
{\twelvepoint\bf \bls1.2\bls
Combinatorial Identities Related to Representations of $\Uqh$
\par}
\vsk1.5>
\VT/
\vsk1.5>
{\it \home/
\vsk.1>
\homeaddr/}
\vsk1.5>
{\sl October \,16, 1998}
\endcenter
\ftext{\=\bls12pt{\tenpoint\sl E-mail\/{\rm:}\homemail/}}
\fi
\ifMag\vsk.4>\else\vsk>\fi
\vsk0>

\sect{Introduction}
Recently N\&Jing discovered in \Cite{J} the following combinatorial identity:
\ifMag\else\vv.3>\fi
$$
\ifMag
\Rline{\sum_{k=0}^\ell
\ \prod_{s=0}^{k-1}\,{\eta^{\)\ell}-\eta^s\over1-\eta^{s+1}\}}\;
\susi\]\Bigl(\]\prod_{1\le a\le k}(t_{\si_a}\!-1)\,
\prod_{k<b\le\ell}(t_{\si_b}\!-\eta^{\)\ell-1})\!\!
\prab{t_{\si_a}\!-\eta\)t_{\si_b}\over t_{\si_a}\!-t_{\si_b}}\Bigr)\)=\,0\,.\!}
\else
\sum_{k=0}^\ell
\ \prod_{s=0}^{k-1}\,{\eta^{\)\ell}-\eta^s\over1-\eta^{s+1}\}}\;
\susi\]\Bigl(\]\prod_{1\le a\le k}(t_{\si_a}\!-1)\,
\prod_{k<b\le\ell}(t_{\si_b}\!-\eta^{\)\ell-1})\!\!
\prab{t_{\si_a}\!-\eta\)t_{\si_b}\over t_{\si_a}\!-t_{\si_b}}\Bigr)\)=\,0\,.
\fi
\global\Tagdef{$\}\II\}$}{Ji}
\Tag{Ji}
$$
\ifMag\else\vvv.3>\fi
In his paper the identity comes from validity of the Serre relations in
some vertex \rep/s of quantum Kac-Moody algebras.
\par
In this note we are going to generalize this identity, see Theorems \[id],
\[idp]. The obtained identities are \eqv/ to existence of a singular vector
in certain tensor products of evaluation \rep/s of the quantum loop algebra
$\Uqh$, see Proposition \[BC], or more generally, of the \eqg/ $\Eqg$.
\par
There are two directions for generalizing the identity \(Ji). First one can
replace linear \fn/s of $\tell$ by \pol/s of larger degree. This is done
in Section \[:Pi], see \(id1).
\par
Furthermore, it is possible to take elliptic theta-\fn/s instead of \pol/s,
see Section \[:Ei]. The resulting identities depend on two extra parameters:
the elliptic modulus $p$ and the dynamical parameter $\al$. In the limit
$p\to 0$ the elliptic \fn/s degenerate into \pol/s and we get a family of
\pol/ identities depending on $\al$ and turning into the identities \(id1)
if either $\al\to 0$ or $\al\to\8$.
\par
The author is grateful to the Max-Planck-Institut f\"ur Mathematik in Bonn,
where the basic part of the paper had been written, for hospitality.

\Sno 0
\Sect[Pi]{Polynomial identities}
Given nonnegative integers $\ell$ and $n$ let $\Pln$ be a set of partitions
\vv.1>
${\la=(\la_1\lc\la_\ell)}$ \st/ ${n\ge{}}\alb{\la_1\lsym\ge\la_\ell\ge1}$.
For a partition $\la$ let $\om_{k,\la}=\#\)\lb\)j\vert\la_j=k\)\rb$.
\vsk.2>
Introduce \indet/s $\tell$, $\xn$, $\yn$, $\eta$. In the paper we use
the following compact notations:
\vv.1>
$$
t=(\tell)\,,\qqq x=(\xn)\,, \qqq y=(\yn)\,.
$$
\vsk.2>
For any $\mn$ set
$$
\gather
\\
\ifMag\cnn-1.3>\else\cnn-1.4>\fi
X_m(u;x;y)\,=\,u\prjm\!(u-y_j)\}\prod_{m<k\le n}\!(u-x_k)\,,
\\
\ifMag\nn7>\else\nn8>\fi
X'_m(u;x;y)\,=\prjm\!(u-x_j)\}\prod_{m<k\le n}\!(u-y_k)\,,
\endgather
$$
and for any \)$\la\in\Pln$ set
$$
\gather
\\
\ifMag\cnn-2.2>\else\cnn-1.9>\fi
r_\la(\eta)\,=\,\pron\,\prod_{s=1}^{\om_{m,\la}}\,{1-\eta\over1-\eta^s}\;,
\\
\nn8>
P_\la(t;x;y;\eta)\,=\,
r_\la(\eta)\)\susi\]\Bigl({}\pral\>X_{\la_a}(t_{\si_a};x;y)\!
\prab{t_{\si_a}\!-\eta\)t_{\si_b}\over t_{\si_a}\!-t_{\si_b}}\Bigr)\,,
\Tag{P}
\\
\ald
\nn8>
P'_\la(t;x;y;\eta)\,=\,
r_\la(\eta)\)\susi\]\Bigl({}\pral\>X'_{\la_a}(t_{\si_a};x;y)\!
\prab{\eta\)t_{\si_a}\!-t_{\si_b}\over t_{\si_a}\!-t_{\si_b}}\Bigr)\,.
\Tag{P'}
\endgather
$$
\ifMag\else\vvv-.2>\fi
Notice that
$$
\gather
\\
\cnn-2.1>
X_m(u;x;y)\,=\,u\>X'_m(u;y;x)
\\
\nn1>
\Text{and}
\nn-6>
P_\la(t;x;y;\eta)\,=\,
\eta^{\)\ell(\ell-1)/2\,-\!\sun\!\om_{m,\la}(\om_{m,\la}-1)/2}\>
t_1\ldots t_\ell\,P'_\la(t;y;x;\eta\1)\,.
\\
\cnn.2>
\endgather
$$
Notice also that both $P_\la$ and $P'_\la$ are \pol/s in all the \indet/s
involved.
\Th{id}
Let $\,x_j=\eta^{\)\ell-1}y_i\,$ for some $\)i<j$. Then
$$
\gather
\sum_{\,\la\in\Plij\!}
c_\la\"{i,j}(x;y;\eta)\,P_\la(t;x;y;\eta)\,=\,0
\Tag{id1}
\\
\nn6>
\Text{where $\;\Plij\]=\lb\>\la=(\la_1\lc\la_\ell)\vert
j\ge\la_1\lsym\ge\la_\ell\ge i\>\rb\,$ and}
\ifMag\nn6>\else\nn2>\fi
{\align
c_\la\"{i,j}(x;y;\eta)\,=\,
(-1)\vpb{\om_{i,\la}}\,\eta\vpb{\,\om_{j,\la}(\om_{j,\la}-1)/2}
\prod_{i<k<j}\!\!\]\prod_{s=0}^{\om_{k,\la}-1}(x_k-\eta^s y_k)\,\x{}&
\ifMag\kern-1.6em\fi
\\
\nn4>
{}\x\,\pral\]\Bigl({}\prod_{i<k<\la_a\!}\!(\eta^{\)\ell-a}\)y_i-x_k)
\prod_{\la_a\]<m<j}(\eta^{\)\ell-a}\)y_i-y_m)\Bigr)\,& .
\ifMag\kern-1.6em\fi
\endalign}
\endgather
$$
\endpro
\nt
The theorem is proved in the next section.
\Ex
Identity \(id1) for $i=1$, $j=n=2$ is \eqv/ to identity \(Ji).
\enddemo
\Rem
Multiplied by a certain \pol/ in $x,y,\eta$ the identity \(id1) becomes more
transparent and understandable, see \(id2).
\mmgood
\enddemo
\Rem
All over the paper we assume integers $\ell$ and $n$ to be fixed. There are
two natural embeddings of $\Pc_{\ell,n-1}$ into $\Pln$ given by either
\>$\la\map\la$ \,or $\>\la\map\la'=(\la_1+1\lc\la_\ell+1)$. Indicating for
a while dependence of \pol/s defined by \(P) on $n$ explicitly, that is,
writing $P_{n,\la}$ instead of $P_\la$, we have
\ifMag\vv-.2>\fi
$$
\align
P_{n,\la}(t;x;y;\eta)\,&{}=\,P_{n-1,\la}(t;x\"n;y\"n;\eta)\>\pral\,(t_a-x_n)\,,
\\
\nn7>
&{}=\,P_{n-1,\la'}(t;x\"1;y\"1;\eta)\>\pral\,(t_a-y_1)\,,
\endalign
$$
where $\>x\"n=(x_1\lc x_{n-1})$, $\>y\"n=(y_1\lc y_{n-1})$,
\vv.1>
$\>x\"1=(x_2\lc x_n)$, $\>y\"1=(y_2\lc y_n)$. Therefore, the claim of Theorem
\[id] for $i=1$, $j=n$ implies the claim of the theorem for general $i<j$.
\enddemo

\Sect[Pf]{Proof of Theorem \[id]}
In this section we think $\xn$, $\yn$, $\eta$ being complex \var/s rather than
\indet/s and consider them as parameters which \fn/s of $\tell$ can in addition
depend on.
\vsk.2>
Let $f,g$ be \pol/s in $\tell$. We define below their scalar product
$\bra f,g\ket\Sf$. Let
\ifMag\vv.1>\fi
$$
S(\tell)\,=\,\pral\,\pron\,(t_a-x_m)\>(t_a-y_m)
\prod_{\tsize{a,b=1\atop a\ne b}}^\ell\}{t_a-\eta\)t_b\over t_a-t_b}\;.
$$
\ifMag\vvv.1>\fi
If $|\eta\)|>1$ and $|x_m|<1$, $|\)y_m|>1$ for all $\mn$, then we set
\ifMag\vv.5>\else\vvn.3>\fi
$$
\bra f,g\ket\Sf\,=\,{1\over(2\pi i)^\ell\>\ell\)!}\;
\int_{\TT^\ell\]}{f(t)\,g(t)\over S(t)}\;\dtt
\vvgood
$$
\ifMag\vvv.2>\fi
\vvn-.2>
where $\,\topsmash{\dtt=\}\pral dt_a/t_a}\,$ and
$\,{\TT^\ell=\lb\)t\in\Cl\vert\ |t_1|=1\llc|t_\ell|=1\)\rb}$ each circle
oriented counterclockwise. Then one can show that $\bra f,g\ket\Sf$ is
a \raf/ of $\xn$, $\yn$, $\eta$.
\vsk.3>
The \pol/s defined by \(P) and \(P') are biorthogonal \wrt/ the introduced
scalar product.
\Th{PP}
\back\Cite{TV1, Theorem C.9\)}
$\ \bra P'_\la\),P_\mu\ket\Sf\)=\)N_\la\1\>\dl_{\la\)\mu}\,$ where
\ifMag\vv.3>\else\vv.1>\fi
$$
N_\la\,=\,\pron\,\prod_{s=1}^{\om_{m,\la}}\,
{(1-\eta^s)\>(x_m-\eta^{s-1}y_m)\over 1-\eta}\;.
\Tag{N}
$$
\endpro
\Pf.
The proof is based on Lemmas \[Pxy] and \[ResI]. The first equality in \(<>)
implies that $\bra P'_\la\),P_\mu\ket\Sf\)=0$ unless $\la\le\mu$, while
the second equality in \(<>) gives $\bra P'_\la\),P_\mu\ket\Sf\)=0$
unless $\la\ge\mu$ and
$\bra P'_\la\),P_\la\ket\Sf\)={P_\la(\yt\la)}\>{P'_\la(\yt\la)}$.
The rest of the proof is straightforward.
\mmgood
\epf
\Rem
Notations here and in \Cite{TV1} are not always the same though we try to keep
them consistent whenever possible. The identification of labels should be
mentioned: a partition $\la$ here corresponds to a label
$\lg=(\om_{1,\la}\)\lc\om_{n,\la})$ in \Cite{TV1}, the \pol/ $P_\la$ being a
numerator of the \tri/ weight \fn/ $w_\lg$ in \Cite{TV1} up to a simple factor.
\enddemo
For any $\la\in\Pln$ introduce points $\xt\la\,,\,\yt\la\in\Cl\}$ as follows:
$$
\gather
\xt\la\,=\,(\eta^{1-\om_{\)1,\la}}x_1\lc x_1,\)\eta^{\)1-\om_{2,\la}}x_2\lc
x_2,\ \ldots\ ,\)\eta^{\)1-\om_{n,\la}}x_n\lc x_n)\,,
\\
\nn6>
\yt\la\,=\,(\eta^{\)\om_{1,\la}-1}\)y_1\lc y_1,\)\eta^{\)\om_{2,\la}-1}\)y_2\lc
y_2,\ \ldots\ ,\)\eta^{\)\om_{n,\la}-1}\)y_n\lc y_n)\,.
\endgather
$$
\vvv.2>
For partitions $\la\),\mu\in\Pln$ say that $\la\ge\mu$ \>if \)$\la_a\ge\mu_a$
for any $\aell$.
\Lm{Pxy}
\bls1.1\bls
$P_\la(\xt\mu)=0$ \)and $P'_\la(\yt\mu)=0$ unless $\la\ge\mu$.
$P_\la(\yt\mu)=0$ \)and $P'_\la(\xt\mu)=0$ unless $\la\le\mu$.
Besides, only the terms in \(P), \(P') corresponding to the identity \perm/
contribute into the values $P_\la(\xt\la)$, $P_\la(\yt\la)$, $P'_\la(\xt\la)$,
$P'_\la(\yt\la)$.
\endpro
\nt
The proof is straightforward.
\Par
\ifMag\else\nt\fi
For a \fn/ $f(\tell)$ and a point ${\ts\}=(\ts_1\lc\ts_\ell)}$ define
\vvm.2>
a multiple residue $\Res\bigl(f(t)\>\dtt\)\bigr)\vst{t=\ts}\]$%
\ifMag\else\nobreak\hskip0pt minus20pt\alb\fi\ by
$$
\gather
\\
\ifMag\cnn->\else\cnn-1.6>\fi
\Res\bigl(f(t)\>\dtt\)\bigr)\vst{t=\ts}=\,\Res\bigl(\,\ldots\,
\Res\bigl(f(\tell)\>(dt_\ell/t_\ell)\bigr)\vst{t_\ell=\ts_\ell}\ldots\,
(dt_1/t_1)\bigr)\vst{t_1=\ts_1}
\\
\ald
\nn4>
\Text{and set}
\nn1>
\xRes(f)\,=\sula\Res\bigl(f(\tell)\>\dtt\)\bigr)\vst{t=x\>\tright\)\la}\,,
\Tag{xRes}
\\
\nn8>
\yRes(f)\,=\sula\Res\bigl(f(\tell)\>\dtt\)\bigr)\vst{t=y\>\trileft\>\la}\,.
\endgather
$$
\Lm{ResI}
\back\Cite{TV1, Lemma C.8}
Let \pol/s $f,g$ be \st/ their product is a \sym/ \pol/ in $\tell$
of degree less than $2n$ in each of the \indet/s and divisible by
$t_1\ldots t_n$. Then
$$
\bra f,g\ket\Sf\,=\,\xRes(fg/S)\,=\,(-1)^\ell\>\yRes(fg/S)\,.
\Tag{<>}
$$
\ifMag\vsk.2>\else\vsk.3>\fi
\endpro
Let $Q_\la$ be a monomial \sym/ \pol/:
$$
Q_\la(\tell)\,=\;{1\over\om_{1,\la}!\ldots\>\om_{n,\la}!}\,\)
\susi t_{\si_1}^{\la_1}\}\ldots\)t_{\si_\ell}^{\la_\ell}\,.
$$
Introduce transition coefficients $A_{\la\)\mu}$, $\la,\mu\in\Pln$:
\ifMag\vvn.2>\fi
$$
P_\la\,=\sumu\!A_{\la\)\mu}\>Q_\mu\,.
\Tag{PQ}
$$
Then by Theorem \[PP] for any $\nu\in\Pln$
\ifMag\vvn.1>\fi
$$
\sumu\!A_{\la\)\mu}\>\bra P'_\nu,Q_\mu\ket\,=\,N_\la\1\>\dl_{\la\)\nu}\,,
$$
and inverting this relation we get
\ifMag\vv.2>\fi
$$
\sula\!N_\la\>\bra P'_\la,Q_\mu\ket\,A_{\la\)\mu}\,=\,\dl_{\mu\)\nu}\,.
$$
Calculating the scalar products by Lemma \[ResI] we have
\ifMag\vv.3>\fi
$$
\sum_{\,\ka,\la\in\Pln\!}\!\!M_\ka\1\>
Q_\mu(\xt\ka)\>P'_\la(\xt\ka)\>N_\la\>A_{\la\)\mu}\,=\,\dl_{\mu\)\nu}\,,
\Tag{MN}
$$
where
$\>M_\ka\1=\}\Res\bigl(S(t)\vpb{-1}\dtt\bigr)\vst{t=x\>\tright\)\ka}\!$.
Notice that $M_\ka$ is a \raf/ in $\xn$, $\eta$ and a \pol/ in $\yn$.
\vsk.2>
The matrix $\bigl[\)Q_\la(\xt\ka)\bigr]_{\ka,\la\in\Pln}\!$ is invertible,
\vv.1>
see e.g.\ \(detQ), and we denote by $B$ the inverse matrix, which is a \raf/
of $\xn$, $\eta$ and, trivially, does not depend on $\yn$. Finally,
the relation \(MN) is transformed to
\ifMag\vv.4>\fi
$$
\sula\!P'_\la(\xt\ka)\>N_\la\>A_{\la\)\mu}\,=\,
M_\ka\>B_{\ka\)\mu}\,.
$$
\Pf of Theorem \[id].
Fix $i<j$. Let $\ka\"j\]=(j\lc j)$. Assume that $y_i=\eta^{\)1-\ell}x_j$.
Then $M_{\ka\"j}=0$, and taking into account the definition \(PQ) of
$A_{\la\)\mu}$ we have an identity
\ifMag\vv-.5>\fi
$$
\sula^{\vp.}\!P'_\la(\xt\ka\"j)\>N_\la\>P_\la\,=\,0\,.
\Tag{id2}
$$
\vvv-.3>
Notice that $\xt\ka\"j=(\eta^{\)1-\ell}x_j\lc x_j)$, so only the term
in \(P') corresponding to the identity \perm/ contribute into the value
$P'_\la(\xt\ka\"j)$, and $P'_\la(\xt\ka\"j)=0$ unless
${\la\in\Plij}\!$.
Calculating ${P'_\la(\xt\ka\"j)\>N_\la}$ explicitly and removing all factors
which does not depend on $\la$ we get the identity \(id1).
Theorem \[id] is proved.
\epf
\Rem
Let us mention two determinant formulae though they are not actually used
in the proof:
\ifMag
\vv-.2>
$$
\align
\quad\det\bigl[Q_\la(\xt\mu)\bigr]_{\la,\mu\in\Pln}=\,
\eta\vpb{-n(n+1)/2\)\cdot\!\!\tsize{n+\ell-1\choose n+1}}
\pron\,\Rlap{x_m}\hp{x}{\vp{\Big|}}^{\!\tsize{n+\ell-1\choose n}}\,\x{}\!\] &
\Tag{detQ}
\\
\nn4>
{}\x\,\prsl\)\prod_{1\le j<k\le n}\!(\eta^sx_k-x_j)\vpb{D(n,\ell,s)} &
\\
\cnn.3>
\endalign
$$
\else
\vv.4>
$$
\det\bigl[Q_\la(\xt\mu)\bigr]_{\la,\mu\in\Pln}=\,
\eta\vpb{-n(n+1)/2\)\cdot\!\!\tsize{n+\ell-1\choose n+1}}\,
\pron\,\Rlap{x_m}\hp{x}{\vp{\Big|}}^{\!\tsize{n+\ell-1\choose n}}\!
\prsl\)\prod_{1\le j<k\le n}\!(\eta^sx_k-x_j)\vpb{D(n,\ell,s)}\kern-.3em
\Tag{detQ}
$$
\vvv.2>
\fi
where
\vv-.3>
$\ \dsize D(n,\ell,s)\,=\!\!\sum_{\tsize{\,r\in\Zp\!\atop 2r\le\ell-|s|-1}}^{}
\!\!{n+\ell-|s|-2r-3\choose n-2}$, \ see formula (A.14) in \Cite{TV1}, \ and
\ifMag\vv->\else\vv-.05>\fi
$$
\det\][A_{\la\)\mu}]\vpp{\)\la,\mu\in\Pln}=\,\pros\,
\prod_{1\le j<k\le n}\!(\eta^s y_j-x_k)\vpb{\tsize{n+\ell-s-2\choose n-1}},
\Tag{detA}
$$
see formula (A.10) in \Cite{TV1}. Formula \(detQ) is a deformation of
a \sym/ power of the Vandermonde determinant. Formula \(detA) implies that
the \pol/s $P_\la$, $\la\in\Pln$ are linear independent if $x,\)y,\)\eta$ are
generic, and there are linear relations between them if \)$x_k=\eta^{\)s}y_j$
\)for some $j<k$ and $s\in\lb\)0\lc\ell-1\)\rb$. For $s=\ell-1$ there is just
one linear relation given by \(id2). Relations for $s<\ell-1$ can be obtained
in a similar manner. They also can be derived from the relation \(id2) written
for \pol/s in less number of \indet/s.
\enddemo

\Sect[Uq]{Tensor products of \emod/s over $\Uqh$}
Let $q$ be a nonzero complex number which is not a root of unity.
Consider the \qg/ $\Uq$ with generators $E,\>F,\>q^H$ and relations
\vvn-.3>
$$
q^H\]E=q\)E\>q^H\,,\qqq q^H\]F=\q F\>q^H\,,\qqq
[\)E\),\]F\>]\;=\,{q^{\)2H}-q^{-2H}\over q-\q}\;,
$$
and the \qlo/ $\Uqh$ with generators $L_{ij}\"{+0},\ L_{ji}\"{-0}\!$,
${1\le j\le i\le 2}$, and $L_{ij}\"s\!$, $i,j=1,2$, $s=\pm 1,\pm 2,\ldots$,
subject to relations \(Ughr).
\par
Let $e_{ij}$, $i,j=1,2$, be the ${2{\x}2}$ matrix with the only nonzero entry
$1$ at the intersection of the \^{$i$-}th row and \^{$j$-}th column. Set
$$
\align
R(u)\,=\,{} &(u\)q-\q)\>(e_{11}\ox e_{11}+e_{22}\ox e_{22})\;+
\\
\ifMag\nn7>\else\nn6>\fi
{} +\,{} &(u-1)\>(e_{11}\ox e_{22}+e_{22}\ox e_{11})\,+\,
u\)(q-\q)\>e_{12}\ox e_{21}\,+\,(q-\q)\>e_{21}\ox e_{12}\,.
\endalign
$$
Introduce the generating series
$L^{\pm}_{ij}(u)=L_{ij}\"{\pm0}+\sum_{s=1}^\8 L_{ij}\"{\pm s}u^{\pm s}\!$.
The relations in $\Uqh$ have the form:
$$
\gather
\\
\ifMag\cnn-1.5>\else\cnn-1.3>\fi
{\align
& R(u/z)\>L^\nu\'1(u)\>L^\nu\'2(z)\,=\,
L^\nu\'2(z)\>L^\nu\'1(u)\>R(u/z)\,,\Rlap{\qqq\nu=\pm\,,}
\Tag{Ughr}
\\
\ifMag\nn8>\else\nn7>\fi
& R(u/z)\>L^+\'1(u)\>L^-\'2(z)\,=\,
L^-\'2(z)\>L^+\'1(u)\>R(u/z)\,,
\endalign}
\\
\ifMag\nn8>\else\nn7>\fi
L_{ii}\"{+0}L_{ii}\"{-0}\,=\,1\,,\qquad L_{ii}\"{-0}L_{ii}\"{+0}\,=\,1\,,
\Rlap{\qqq i=1,2\,,}
\endgather
$$
where \,$L^\nu\'1(u)=\sum_{ij}\)e_{ij}\ox 1\ox L_{ij}^\nu(u)$ \,and
\,$L^\nu\'2(u)=\sum_{ij}\)1\ox e_{ij}\ox L_{ij}^\nu(u)$.
\vsk.15>
The \qlo/ $\Uqh$ is a Hopf algebra with a coproduct
$$
\Dl:L^\pm_{ij}(u)\,\map\,\tsum_k\)L^\pm_{ik}(u)\ox L^\pm_{kj}(u)\,.
$$
\ifMag\vv-.8>\else\vv-.3>\fi
There is a one-parametric family of \aut/s
\mmgood
\>$\rho_z:L_{ij}^\pm(u)\,\map\,L_{ij}^\pm(u/z)$ \,and the evaluation \hom/
$\epe:\Uqh\to\Uq$:
\ifMag\else\vv-.6>\fi
$$
\epe\):\)L^\pm(u)\,\map\,{}\mp\,u^{\thi_\pm}
\pmatrix
\,u\)q^H-\qH\; & u\)F\>(q-\q)\> \\ \nn8>
\)E\>(q-\q)\; & u\)\qH-q^H
\endpmatrix\kern-1em
$$
where ${\)\thi_+=0}$, ${\)\thi_-=-1\)}$. For any \Umod/ $V$ denote by $V(z)$
the \Uhmod/ obtained from $V$ via the \hom/ ${\epe\o\rho_z}$. The module $V(z)$
is called the \emod/. In a tensor product $\Voz$ of \emod/s actions of
the series $L^+(u)$ and $(-u)^nL^-(u)$ coincide, so it is enough to look
at the action of only one of them. Besides, this implies that the series
$L^+(u)$ acts as a \pol/ in $u$, so it can be evaluated for $u$ being
a complex number.
\vsk.2>
Let $V_1\lc V_n$ be $\Uq$ \Vmod/s with \hw/s $q^{\)\La_1}\lc q^{\)\La_n}$ and
\gv/s $v_1\lc v_n$, \resp/. Using commutation relations in $\Uqh$ it is rather
straightforward to obtain the following statements.
\Prop{BC}
Consider $\Uqh$ \emod/s $V_1(z_1)\lc V_n(z_n)$.
\vv.1>
Let $z_i=q^{\)2\La_i+2\La_j-2\ell}z_j$ for some $i<j$ and $\ell\in\Zp$. Then
\atem The vector ${\vox\in\Voz}$ generates
\ifMag\vv.25>\else\vv.1>\fi
a \$\Uqh$-submodule, which is annihilated by the action of a product
\>$L^+_{21}(q^{\)2\La_j}z_j)\ldots L^+_{21}(q^{\)2\La_j-2\ell}z_j)$.
\vsk.3>
\bitem
\ifMag\vv.15>\else\vv.1>\fi
A vector ${\vti\>=L^+_{12}(q^{\)2\La_j}z_j)\ldots
L^+_{12}(q^{\)2\La_j-2\ell}z_j)\,\voxx\in\Vozz}$
is singular \wrt/ the \$\Uqh$-action, that is, \>$L^\pm_{21}(u)\,\vti=0$.
\vvgood
\endpro
\Cr{BCC}
Let $z_i=q^{\)2\La_i+2\La_j-2\ell}z_j$ for some $i<j$ and $\ell\in\Zp$. Then
for arbitrary $\tell$
$$
L^+_{21}(q^{\)2\La_j}z_j)\ldots L^+_{21}(q^{\)2\La_j-2\ell}z_j)\>
L^+_{12}(t_1)\ldots L^+_{12}(t_\ell)\,\vox\,=\,0
\Tag{BC1}
$$
in $\Voz$ and
$$
L^+_{21}(t_1)\ldots L^+_{21}(t_\ell)\>L^+_{21}(q^{\)2\La_j}z_j)\ldots
L^+_{12}(q^{\)2\La_j-2\ell}z_j)\,\voxx\,=\,0
\Tag{BC2}
$$
in $\Vozz$.
\endpro
\nt
Both relations \(BC1) and \(BC2) are \eqv/ to the identitiy \(id2)
and, hence, the identity \(id1). For \(BC1) it follows from Proposition \[KBI]
and formula \(N), while for \(BC2) an analogue of Proposition \[KBI] for
the tensor product $\Vozz$ is required.
\Rem
It is rather easy to see that Corollary \[BCC], Proposition \[KBI] and
the determinant formula \(detA) together imply Proposition \[BC].
\enddemo
Let parameters $q$, $q^{\)\La_1}\lc q^{\)\La_n}$, $\zn$ used in this section
be related to the previously used parameters $\eta$, $\xn$, $\yn$ as follows:
$$
\eta=q^2\,,\qqq x_m=\)q^{\)2\La_m}z_m\,,\qqq y_m=\)q^{-2\La_m}z_m\,.
\kern-2em
\mmgood
$$
\Prop{KBI}
\back\Cite{KBI},\;\Cite{TV1, Lemma 4.3},\;\Cite{TV2, Lemma 4.18}
\ifMag
$$
\gather
\\
\cnn-1.2>
{\align
L^+_{12}(t_1)\ldots L^+_{12}(t_\ell)\,\vox\,=\,
(q &{}-\q)^\ell\>\tpron(-z_m)^{-\ell}\!\sula\!P_\la(\tell)\,\x{}
\\
\nn6>
{}\x\,\prod_{1\le j<k\le n}\! &
q^{\)\La_j\om_{k,\la}-\)\La_k\om_{j,\la}-\>\om_{j,\la}\om_{k,\la}}\,\Fv\,.
\endalign}
\\
\nn15>
\Lline{L^+_{21}(t_1)\ldots L^+_{21}(t_\ell)\,\Fv\,={}}
\\
\nn8>
\Rline{{}=\tpron(-z_m)^{\)\om_{m,\la}-\)\ell}\,\)
P'_{\la}(\tell)\,\pron\,\prod_{s=1}^{\om_{m,\la}}\,
{(q^s-q^{-s})\>(q^{\)2\La_m-s+1}-q^{-2\La_m+s-1})\over q-\q}\;\x{}\!}
\\
\nn6>
\Rline{{}\x\,\prod_{1\le j<k\le n}\!
q^{\)\La_k\om_{j,\la}-\)\La_j\om_{k,\la}-\>\om_{j,\la}\om_{k,\la}}\,\vox\,.\!}
\\
\cnn->
\endgather
$$
\else
$$
\gather
{\align
L^+_{12}(t_1)\ldots L^+_{12}(t_\ell)\,\vox\, &{}=\,
(q-\q)^\ell\>\tpron(-z_m)^{-\ell}\sula\!P_\la(\tell)\,\x{}
\\
\nn6>
& {}\>\x\,\prod_{1\le j<k\le n}\!
q^{\)\La_j\om_{k,\la}-\)\La_k\om_{j,\la}-\>\om_{j,\la}\om_{k,\la}}\,\Fv\,.
\endalign}
\\
\nn12>
\alignedat2
L^+_{21}(t_1)\ldots L^+_{21} & (t_\ell)\,\Fv\,={} &&
\\
\nn7>
& {}=\tpron(-z_m)^{\)\om_{m,\la}-\)\ell}\,
P'_{\la}(\tell)\,\pron\,\prod_{s=1}^{\om_{m,\la}}\,
{(q^s-q^{-s})\>(q^{\)2\La_m-s+1}-q^{-2\La_m+s-1})\over q-\q}\;\x{} &&
\\
\nn6>
&& \Llap{{}\x\,\prod_{1\le j<k\le n}\!
q^{\)\La_k\om_{j,\la}-\)\La_j\om_{k,\la}-\>\om_{j,\la}\om_{k,\la}}\,\vox\,.} &
\endalignedat
\\
\cnn->
\endgather
$$
\fi
\endpro

\Sect[Ei]{Elliptic identities}
In this section we extend Theorem \[id] to the elliptic case, see Theorem
\[idp]. The obtained identities have the same relation to \rep/s of the \eqg/
$\Eqg$ as the identities \(id1) have to \rep/s of the \qlo/ $\Uqh$.
\Par
Fix a nonzero complex number $p$ \st/ ${|p|<1}$. Let
\ifMag
\vv.1>
$$
(u)\9\>=\,(u;p)\9\>=\,\topsmash{\tprod_{s=0}^\8(1-p^su)}
$$
\vvv-.3>
\else
$(u)\9\]=(u;p)\9\]=\topsmash{\prod_{s=0}^\8(1-p^su)}$
\vv-.4>
\fi
and let ${\tht(u)=\tht(u;p)=(u)\9(p\)u\1)\9(p)\9}$ be the Jacobi theta-\fn/.
\vsk.3>
In addition to $\tell$, $\xn$, $\yn$ and $\eta$ introduce one more \var/
$\al$ called the \em{dynamical parameter}.
\Rem
The parameter $\al$ is related to the dynamical \var/ in the \eqg/ $\Eqg$.
\enddemo
For any $\mn$ set
$$
\gather
\\
\cnn-1.4>
Z_m(u;x;y;\al)\,=\,\tht\bigl(\al_m\1 u/x_m)\}
\prod_{1\le j<m}\!\tht(u/y_j)\}\prod_{m<k\le n}\!\tht(u/x_k)\,,
\\
\nn9>
Z'_m(u;x;y;\al)\,=\,\tht\bigl(\al_m\)u/y_m)\}
\prod_{1\le j<m}\!\tht(u/x_j)\}\prod_{m<k\le n}\!\tht(u/y_k)\,,
\endgather
$$
\ifMag\vvv.2>\fi
where $\>\smash{\al_m=\al\tprjm x_j/y_j}$, \;and for any \)$\la\in\Pln$ set
$$
\gather
\\
\ifMag\cnn-1.2>\else\cnn-1.3>\fi
\rho_\la(\eta)\,=\,\pron\,
\prod_{s=1}^{\om_{m,\la}}\,{\tht(\eta)\over\tht(\eta^s)}\;,
\\
\ifMag\ald\nn7>\else\nn8>\fi
\Xj_\la(t;x;y;\eta\);\al)\,=\,\rho_\la(\eta)\)\susi\]\Bigl({}
\pral\>Z_{\la_a}(t_{\si_a};x;y;\al\)\eta^{\)2a-2\ell})\!\prab\!
{\tht(\eta\)t_{\si_b}/t_{\si_a})\over\tht(t_{\si_b}/t_{\si_a})\!}\,\Bigr)\,,
\ifMag\kern-3em\fi
\Tag{T}
\\
\nn8>
\Xj'_\la(t;x;y;\eta\);\al)\,=\,\rho_\la(\eta)\)\susi\]\Bigl({}
\pral\>Z'_{\la_a}(t_{\si_a};x;y;\al\)\eta^{\)2\ell-2a})\!\prab\!
{\tht(\eta\)t_{\si_a}/t_{\si_b})\over\tht(t_{\si_a}/t_{\si_b})\!}\,\Bigr)\,.
\ifMag\kern-3em\fi
\Tag{T'}
\endgather
$$
\vvv-.4>
Notice that
$$
\gather
\\
\cnn-2.3>
Z_m(u;x;y;\al)\,=\,Z'_m(u;y;x;\al\1)
\\
\nn2>
\Text{and}
\nn-8>
\Xj_\la(t;x;y;\eta\);\al)\,=\,
\eta^{\)\ell(\ell-1)/2\,-\!\sun\!\om_{m,\la}(\om_{m,\la}-1)/2}\>
\Xj'_\la(t;y;x;\eta\1;\al\1)\,.
\endgather
$$
\Th{idp}
Let $\,x_j=\eta^{\)\ell-1}y_i\,$ for some $\)i<j$. Then
\ifMag
$$
\gather
\sum_{\,\la\in\Plij\!}
C_\la\"{i,j}(x;y;\eta\);\al)\,\Xj_\la(t;x;y;\eta\);\al)\,=\,0
\Tag{idp1}
\\
\nn6>
\Text{where $\;\Plij\]=\lb\>\la=(\la_1\lc\la_\ell)\vert
j\ge\la_1\lsym\ge\la_\ell\ge i\>\rb\,$,}
\nn9>
\alignedat2
C_\la\"{i,j} & \Rlap{(x;y;\eta\);\al)\,={}} &&
\\
\nn1>
{}=\,{}&&
(\al_i\)x_i/y_i)\vpb{\om_{i,\la}}\,\eta\vpb{-\om_{i,\la}(\om_{i,\la}-1)}
\prod_{i<k<j}\!\!\]\prod_{s=0}^{\om_{k,\la}-1}{\tht(\eta^{-s}x_k/y_k)\over
\tht(\eta^s\al_{k,\la}\1)\>\tht(\eta^{1-s-\om_{k,\la}}\al_{k,\la}\>x_k/y_k)}\;
\x{}&
\\
\nn7>
&&\Llap{{}\x{}}\prod_{s=0}^{\om_{i,\la}-1}
{1\over\tht(\eta^{1-s-\om_{i,\la}}\al_{i,\la}\>x_i/y_i)}
\prod_{s=0}^{\om_{j,\la}-1}{1\over\tht(\eta^s\al_{j,\la}\1)}\;
\prod_{\,a=\om_{j,\la}+1\!}^{\,\ell-\om_{i,\la}\!}\!
\tht(\al_{\la_a}\:\eta^{a-\ell}\)y_i/y_{\la_a}\:)\;\x{}&
\\
\nn5>
&&\pral\]\Bigl({}\prod_{i<k<\la_a\!}\!\tht(\eta^{\)\ell-a}\)y_i/x_k)
\prod_{\la_a\]<m<j}\tht(\eta^{\)\ell-a}\)y_i/y_m)\Bigr)\,&,
\endalignedat
\endgather
$$
\else
$$
\gather
\sum_{\,\la\in\Plij\!}
C_\la\"{i,j}(x;y;\eta\);\al)\,\Xj_\la(t;x;y;\eta\);\al)\,=\,0
\Tag{idp1}
\\
\nn7>
\Text{where $\;\Plij\]=\lb\>\la=(\la_1\lc\la_\ell)\vert
j\ge\la_1\lsym\ge\la_\ell\ge i\>\rb\,$,}
\nn8>
{\align
C_\la\"{i,j}(x;y;\eta\);\al)\,=\,
(\al_i\)x_i/y_i)\vpb{\om_{i,\la}}\,\eta\vpb{-\om_{i,\la}(\om_{i,\la}-1)}
\prod_{i<k<j}\!\!\]\prod_{s=0}^{\om_{k,\la}-1}{\tht(\eta^{-s}x_k/y_k)\over
\tht(\eta^s\al_{k,\la}\1)\>\tht(\eta^{1-s-\om_{k,\la}}\al_{k,\la}\>x_k/y_k)}
\,\x{}&
\\
\nn6>
{}\x\prod_{s=0}^{\om_{i,\la}-1}
{1\over\tht(\eta^{1-s-\om_{i,\la}}\al_{i,\la}\>x_i/y_i)}
\prod_{s=0}^{\om_{j,\la}-1}{1\over\tht(\eta^s\al_{j,\la}\1)}
\prod_{\,a=\om_{j,\la}+1\!}^{\,\ell-\om_{i,\la}\!}\!
\tht(\al_{\la_a}\:\eta^{a-\ell}\)y_i/y_{\la_a}\:)\;\x{}&
\\
\nn5>
{}\x\,\pral\]\Bigl({}\prod_{i<k<\la_a\!}\!\tht(\eta^{\)\ell-a}\)y_i/x_k)
\prod_{\la_a\]<m<j}\tht(\eta^{\)\ell-a}\)y_i/y_m)\Bigr)\,& ,
\endalign}
\\
\cnn-.2>
\endgather
$$
\fi
$\al_{k,\la}=\al\!\prod_{1\le j<k}^{\vp.}\!\eta^{-2\)\om_{j,\la}}\)x_j/y_j$
\,\)and \;$\al_k=\al\!\prod_{1\le j<k}\!x_j/y_j$.
\ifMag\vv.3>\else\vv.5>\fi
\endpro
\Ex
Identity \(idp1) for $i=1$, $j=n=2$ takes the form
$$
\alignat2
\\
\ifMag\cnn-1.3>\else\cnn-1.4>\fi
& \sum_{k=0}^\ell\,(-1)^k\>\tht(\eta^{\)2k}\bt)\,
\prod_{s=0}^{k-1}\,{\eta^s\>\tht(\eta^{\)\ell-s})\>\tht(\eta^s\bt)\over
\tht(\eta^{s+1})\>\tht(\eta^{s+\ell+1}\bt)}\;\x{} &&
\ifMag\kern-1.5em\else\kern-.5em\fi
\Tag{idp2}
\\
\nn12>
&{}\x\susi\]\Bigl(\]\prod_{1\le a\le k}\tht(t_{\si_a})\>
\tht(\eta^{\)2-2a-\ell}t_{\si_a}/\]\bt)\;\prod_{k<b\le\ell}
\tht(\eta^{\)1-\ell}t_{\si_b})\>\tht(\eta^{\)1-2b}t_{\si_b}/\]\bt)\;\x{} &&
\ifMag\kern-1.5em\else\kern-.5em\fi
\\
\nn1>
&&\Llap{{}\x\prab
{\tht(\eta\)t_{\si_b}/t_{\si_a})\over\tht(t_{\si_b}/t_{\si_a})}
\Bigr)\,}{}\}=\,0\,.\ifMag\kern-1.5em\else\kern-.5em\fi
\\
\ifMag\noalign{\goodbm}\fi
\cnn-.1>
\endalignat
$$
Here $\>\bt=\eta^{\)1-2\ell}\al\>x_1/y_1$. If $p\to 0$ and then either
$\bt\to 0$ or $\bt\to\8$, then \(idp2) transforms into \(Ji).
\enddemo
\Pf of Theorem \[idp].
The proof of Theorem \[idp] is very similar to the proof of Theorem \[id].
We descripe the main steps below. We do not indicate dependence on $\xn$,
$\yn$, $\eta$, $\al$ explicitly.
\goodbreak
\Par
Denote by $\EE$ the space of \fn/s $f(\tell)$ \hol/ outside the coordinate
hyperplanes $t_a=0$, $\aell$, and \st/
\ifMag\vv-.2>\fi
$$
\gather
f(t_1\lc p\)t_a\lc t_\ell)\,=\,
t_a^{-2n}\!\tpron(x_m\>y_m)\,f(\tell)\,.
\\
\nn1>
\Text{Let}
\nn-6>
\Om(\tell)\,=\,\pral\,\pron\,\tht(t_a/x_m)\>\tht(t_a/y_m)
\prod_{\tsize{a,b=1\atop a\ne b}}^\ell
{\tht(\eta\)t_a/t_b)\over\tht(t_a/t_b)}\;.
\\
\cnn-.5>
\endgather
$$
There is an analogue of Lemma \[ResI].
\Lm{Res}
\back\Cite{TV1, Lemma C.11}
For any \fn/ $f\in\EE$ we have
$$
\xRes(f/\Om)\,=\,(-1)^\ell\>\yRes(f/\Om)\,.
$$
\endpro
Define a scalar product $\bra f,g\ket\Of=\xRes(fg/\Om)$.
The \fn/s \(T) and \(T') are biorthogonal \wrt/ the introduced scalar product.
\Th{XX}
\back\Cite{TV1, Theorem C.9\)}
$\ \bra \Xj'_\la\),\Xj_\mu\ket\Of\)=\)D_\la\1\>\dl_{\la\)\mu}\,$ where
\ifMag\vv.1>\fi
$$
D_\la\,=\,(-1)^\ell\pron\!\prod_{s=0}^{\om_{m,\la}-1}
{(p)^3_\8\>\tht(\eta^{s+1})\>\tht(\eta^{-s}x_m/y_m)\over\tht(\eta)\>
\tht(\eta^s\al_{m,\la}\1)\>\tht(\eta^{1-s-\om_{m,\la}}\al_{m,\la}\>x_m/y_m)}
\Tag{D}
$$
and $\,{\al_{m,\la}=\al\!\prjm^{\vp.}\!\eta^{-2\)\om_{\la,j}}\)x_j/y_j}$.
\ifMag\vv>\fi
\endpro
Let $\thi_1\lc \thi_n$ be the following \fn/s in one \var/:
$$
\thi_m(u)\,=\,u^{m-1}\>\tht\bigl(-\)p^{\)m-1}\)\eta^{\)\ell-1}\)\al\1\!\]
\tpron\!x_m\>(-u)^n\);\)p^n\bigr)\,(p^n;p^n)\1_\8\,(p;p)^n_\8\,.
\Tag{thi}
$$
They are linearly independent, since $\thi_m(\eps u)=\eps^{m-1}\>\thi(u)$
where $\eps=e^{2\pi i/n}\!$. For any $\la\in\Pln$ set
$$
\Tho_\la(\tell)\,=\;{1\over\om_{1,\la}!\ldots\>\om_{n,\la}!}\,\)
\susi\thi_{\la_1}(t_{\si_1})\ldots\)\thi_{\la_n}(t_{\si_\ell})\,.
\Tag{Tho}
$$
The \fn/s $\Tho_\la$, $\la\in\Pln$, are linearly independent, and
the \fn/s $\Xj_\mu$ are linear combinations of the \fn/s $\Tho_\la$,
see Proposition \[XT].
\goodbm
\par
We proceed further like in Section \[:Pf] and under the assumption of
Theorem \[idp]\): $y_i=\eta^{\)1-\ell}x_j$ for some $i<j$, we get an identity
$$
\sula^{\vp.}\!\Xj'_\la(\xt\ka\"j)\>D_\la\>\Xj_\la\,=\,0\,,
$$
and then the identity \(idp1). Theorem \[idp] is proved.
\epf
\Rem
There are elliptic analogues of determinant formulae \(detQ) and \(detA),
see formulae (B.11) and (B.8) in \Cite{TV1}\):
$$
\align
\\
\ifMag\cnn-.9>\else\cnn-1.4>\fi
\det\bigl[\Tho_\la(\xt\mu)\bigr]_{\la,\mu\in\Pln} &{}=\,K\>
\eta\vpb{n(1-n)/2\)\cdot\!\!\tsize{n+\ell-1\choose n+1}}\,
\pron\,(-x_m)\vpb{(m-1)\tsize{n+\ell-1\choose n}}\,\x{} \kern-1em
\Tagg{detT}
\\
\ifMag\nn9>\else\nn8>\fi
& {}\>\x\,\pros\,\tht(\eta^s\)\al\1)\vpb{\tsize{n+s-1\choose n-1}}
\}\prsl\)\prod_{1\le j<k\le n}\!\tht(\eta^sx_j/x_k)\vpb{D(n,\ell,s)}
\ifMag\kern-2.5em\fi
\\
\ifMag\cnn.5>\else\cnn.3>\fi
\endalign
$$
where \ $\dsize K\,=\,\Bigl[\>\Rph{(p)}{(p)\9}\vpb{n^2-1}\>\prmn\>
\Bigl(\>{\tht(e^{2\pi i\)m/n})\over e^{2\pi i\)m/n}\}-1}\)\Bigr)
^{\!\!\raise 3pt\mbox{\ssize m-n}}\>\Bigr]^{\tsize{n+\ell-1\choose n}}\!$,\quad
and
$$
\align
\det\][\Ae_{\la\)\mu}]\vpp{\)\la,\mu\in\Pln}& {}=\,K\1\}\prsl\)\prmn\>
\tht({\tsize\eta^{s+\ell-1}\al\1\!\!\prod_{1\le j\le m}\!y_j/x_j})
\vpb{d(n,m,\ell,s)}\,\x{} \kern-1em
\Tagg{detAe}
\\
\nn10>
& \>{}\x\,\pron\,y\_m\vpb{(m-n){\tsize{n+\ell-1\choose n}}}\;\pros\)
\prod_{1\le j<k\le n}\!\tht(\eta^s\)y_j/x_k)
\vpb{\tsize{n+\ell-s-2\choose n-1}} \ifMag\kern-1.5em\fi
\\
\cnn.5>
\endalign
$$
where \ $\dsize\Xj_\la\,=\sumu\!\Ae_{\la\)\mu}\>\Tho_\mu$ \ and
\ $\dsize d(n,m,\ell,s)\,=
\sum_{\tsize{\vp( i,j\ge 0\atop\tsize{\vp( i+j<\ell\atop i-j=s\vp(}}}
{m-1+i\choose m-1}\,{n-m-1+j\choose n-m-1}$.
\enddemo
\Rem
Consider a limit $p\to 0$ in all formulae of this section, which essentially
amounts to replacing the theta-\fn/ $\tht(u)$ by the linear \fn/ $1-u$.
The \fn/s $\thi_1\lc\thi_n$ tend to power \fn/s as $p\to 0$:
$$
\thi_1(u)\to 1+\eta^{\)\ell-1}\)\al\1\!\]\tpron\!x_m\>(-u)^n\,,\qqq
\thi_m(u)\to u^{m-1}\,,\qquad m=2\lc n\,.
$$
The limit ${p\to 0}$ of the identities \(idp1) and the determinant formulae
\(detT) and \(detAe) deliver one-parametric deformations of the identities
\(id1) and the determinant formulae \(detQ) and \(detA), which can be recovered
in the limit either $\al\to 0$ or $\al\to\8$. The limit $\al\to 0$ is
straightforward while the limit $\al\to\8$ produces formulae \(id1), \(detQ),
\(detA) up to a certain change of \var/s.
\enddemo
\Prop{XT}
For any $\la\in\Pln$ the \fn/ $\Xj_\la$, \cf. \(T), is a linear combination
of the \fn/s $\Tho_\mu$, $\mu\in\Pln$, \cf. \(Tho).
\endpro
\Pf.
Let $\Ec$ be a space of \fn/s $f(u)$ \hol/ for $u\ne 0$ and \st/
$$
f(p\)u)\,=\,\al\)\eta^{1-\ell}\]\tpron\}x_m\,(-u)^n\,f(u)\,.
$$
Let $\Ec_\ell$ be a space of \sym/ \fn/s $f(\tell)$ which considered as \fn/s
of one \var/ $t_a$ belong to $\Ec$ for any $\aell$. In particular, $\Ec_1=\Ec$.
\par
It is clear that $\dim\Ec=n$, say by Fourier series. A basis in $\Ec$ is given
by the \fn/s $\thi_1\lc\thi_n$. Hence, the \fn/s $\Tho_\mu$, $\mu\in\Pln$, form
a basis in the space $\Ec_\ell$.
\par
Since for any $\la\in\Pln$ the \fn/ $\Xj_\la$ belongs to the space $\Ec_\ell$,
the proposition follows.
\epf

\myRefs
\widest{VV}
\parskip.1\bls

\ref\Key J
\by N\&Jing
\paper Quantum Kac-Moody algebras and vertex \rep/s
\jour \LMP/ \vol 44 \yr 1998 \pages 261--271 % \issue 2
\endref

\ref\Key KBI
\by \Kor/, N\&M\&Bogolyubov and A\&G\&Izergin
\book Quantum inverse scattering method and correlation \fn/s
\yr 1993 \publ \CUP/ % \pages 1--556
\endref

\ref\Key T
\by \VT/
\paper Bilinear identity for \$q$-\hint/s
\jour Preprint \yr 1997 \pages 1--24 \info to appear in Osaka J. Math.
\endref

\ref\Key TV1
\by \VT/ and \Varch/
\paper Geometry of \$q$-\hgeom/ \fn/s, \qaff/s and \eqg/s
\jour Ast\'erisque \vol 246 \yr 1997 \pages 1--135
\endref

\ref\Key TV2
\by \VT/ and \Varch/
\paper Geometry of \$q$-\hgeom/ \fn/s as a bridge between Yangians and
\qaff/s \jour \Inv/ \yr 1997 \vol 128 \issue 3 \pages 501--588
\endref

\endRefs

\bye